\documentclass[a4paper,11pt]{article}

\usepackage[utf8]{inputenc}
\usepackage{amsmath}
\usepackage{amsthm}
\usepackage{amssymb}
\usepackage{amsfonts}
\usepackage{dsfont}
\usepackage{natbib}
\usepackage{verbatim}
\usepackage{graphicx}
\usepackage[font={stretch=1.0,footnotesize}]{caption}
\usepackage{float}
\usepackage{enumerate}
\usepackage[margin=1in]{geometry}
\usepackage{multicol}
\usepackage{setspace}
\usepackage{hyperref}
\usepackage{tikz}
\usepackage{tikz-qtree}
\usetikzlibrary{bayesnet}
\usepackage{framed}
\def \dsP {\text{$\mathds{P}$}}
\def \dsE {\text{$\mathds{E}$}}

\def \dsN {\text{$\mathds{N}$}}
\def \dsR {\text{$\mathds{R}$}}
\def \dsOne {\text{$\mathds{1}$}}

\theoremstyle{plain}
\newtheorem{theorem}{Theorem}[section]

\newtheorem{proposition}[theorem]{Proposition}
\newtheorem{corollary}[theorem]{Corollary}

\theoremstyle{remark}
\newtheorem{definition}[theorem]{Definition}

\title{Posterior Concentration Rates for Bayesian Penalized Splines}
\author{Paul Bach and Nadja Klein$^{\ast}$}

\begin{document}
\maketitle
\thispagestyle{empty}

\begin{abstract}
   \noindent Despite their widespread use in practice, the asymptotic properties of Bayesian penalized splines have not been investigated so far. We close this gap and study posterior concentration rates for Bayesian penalized splines in a Gaussian nonparametric regression model. A key feature of the approach is the hyperprior on the smoothing variance, which allows for adaptive smoothing in practice but complicates the theoretical analysis considerably as it destroys conjugacy and precludes analytic expressions for the posterior moments. To derive our theoretical results, we rely on several new concepts including a carefully defined proper version of the partially improper penalized splines prior as well as an innovative spline estimator that projects the observations onto the first basis functions of a Demmler-Reinsch basis. Our results show that posterior concentration at near optimal rate can be achieved if the hyperprior on the smoothing variance strikes a fine balance between oversmoothing and undersmoothing, which can for instance be met by a Weibull hyperprior with shape parameter 1/2. We complement our theoretical results with empirical evidence demonstrating the adaptivity of the hyperprior in practice. 
\end{abstract}


\noindent
\vfill
\noindent{\small $^\ast$Paul Bach is a PhD Student at Humboldt-Universit\"at zu Berlin. Nadja Klein is Emmy Noether Research Group Leader in Statistics and Data Science at Humboldt-Universit\"at zu Berlin. Correspondence should be directed to~Nadja Klein, Humboldt-Universit\"at zu Berlin, Unter den Linden 6, 10099 Berlin. Email: nadja.klein@hu-berlin.de
\\
\noindent \textbf{Acknowledgments}: Paul Bach and Nadja Klein gratefully acknowledge support  by the German research foundation (DFG) through the  Emmy Noether grant KL 3037/1-1. }

\setcounter{page}{0}
\pagestyle{plain}

\newpage


\section{Introduction}\label{sec:Intro}

Bayesian penalized splines are extremely popular in applied regression modelling: Since their introduction in the early 2000's by \citet{BerCarRup2002,LanBre2004}, they have been applied in numerous settings ranging from nonparametric and additive regression over generalized additive models, additive survival analysis and multivariate response models to structured additive distributional and quantile regression \citep{LanBre2004,BreLan2006,HenBreFah2006,LanAdeFahSte2003,KleKneLanSoh2015,WalKneYue2013}. Further applications include measurement error models and transformation models \citep{BerCarRup2002,SonLu2012}; see also \citet{FahKneLan2013,HarRupWan2018} for an overview and further real data examples.

The basic idea of the approach is to expand the unknown function to be estimated in a rich spline basis with a large number of knots and to use a prior on the spline coefficients that penalizes the roughness of the associated spline function in order to prevent overfitting. The overall amount of regularization is controlled adaptively through a hyperprior $p(\tau^2)$ on the smoothing variance $\tau^2$. 

Despite its overwhelming success in applied statistics, the theoretical understanding of the approach is lagging far behind. In particular, the asymptotic properties of Bayesian penalized splines have not been investigated so far. In this paper, we make an important step to close this gap in the literature and study posterior concentration rates for Bayesian penalized splines in a Gaussian nonparametric regression model.

The derivation of posterior concentration rates for Bayesian penalized splines is already for the Gaussian response model that we study in this paper quite challenging. This is because of the hyperprior $p(\tau^2)$ on the smoothing variance, which is a key feature of the approach and both a blessing and a curse: From a practical perspective, the hyperprior is a blessing as it allows to adapt to the smoothness of the unknown function. From a theoretical perspective, however, the hyperprior is a curse as it destroys conjugacy and precludes analytic expressions for the posterior moments, which considerably complicates the asymptotic analysis. 

The derivation of posterior concentration rates is generally much more involved if analytic expressions for the posterior moments are unavailable \citep[cf.~][Chapter 8]{Ghovan2017}. In this case one cannot use a ``direct approach'', where one shows that the posterior mean converges to the truth and that the posterior variance vanishes sufficiently fast. To overcome this difficulty, we use the Kullback-Leibler (KL) and testing strategy for the proof of our main posterior concentration result. The KL and testing strategy is a general strategy to prove posterior concentration results that does not require analytic expressions for the posterior moments \citep{Sch1965,Bar1988,GhoGhovan2000,Ghovan2007,Rou2016}. While the KL and testing strategy itself is well-established as a general proof technique, we need to introduce several completely new tools when applying it in the present context of Bayesian penalized splines.

First, the KL and testing strategy requires a proper prior but the penalized splines prior is partially improper. Previous definitions of a proper prior \citep[e.g.][]{WanOrm2008} are based on a mixed model type reparametrization \citep[][]{FahKneLan2004b} but these are not suited for the derivation of posterior concentration rates. Therefore, we introduce a new definition of a proper prior that is based on empirical projections and allows us to use the KL and testing strategy in the present setting.  Second, the KL and testing strategy requires the definition of a suitable test. Given the close relation of Bayesian penalized splines and frequentist penalized splines, it is a natural idea to define this test on the basis of the frequentist penalized splines estimator. The frequentist penalized splines estimator was introduced by \cite{OS1986} and its asymptotic behavior has been studied extensively \citep[e.g.][]{ClaKriOps2009,KauKriFah2009,Xia2019}. However, the test based on the frequentist penalized splines estimator is not suitable for the asymptotic analysis of Bayesian penalized splines because the type II error probability of this test is extremely difficult to control. Therefore, our test relies on a  new spline estimator, which we refer to as \emph{truncated DR estimator}. This estimator projects the observations onto the first basis functions of a \emph{Demmler-Reinsch (DR) basis}, which is a special orthonormal basis of the space of splines. As a byproduct of our posterior concentration results, we establish that the truncated DR estimator is able to achieve the optimal rate of convergence \citep{Sto1982}.

\subsection{Related literature} 

Several papers have derived posterior concentration rates for Bayesian spline regression approaches without roughness penalty \citep[e.g.][]{Ghovan2007,Jonvan2012,YooGho2016,BaiMorAnt2020,WeiReiHop2020}. In Section 7.7 of their seminal paper on posterior concentration rates for non-iid observations, \cite{Ghovan2007} consider a Gaussian nonparametric regression model and expand the unknown function in a relatively small B-spline basis with iid standard Gaussian priors on the coefficients. The work of \cite{Jonvan2012} extends this approach to multivariate functions using tensor product B-splines. A similar approach is used by \cite{YooGho2016}, who additionally consider contraction with respect to the uniform norm. 
\cite{YooGho2016} use a Normal-Inverse-Gamma prior, which allows them to exploit conjugacy for their derivation of posterior concentration rates. 

Due to the hyperprior $p(\tau^2)$ on the smoothing variance, we cannot exploit conjugacy for our derivation of posterior concentration rates and the formal structure of the problem that we face shows parallels to that of shrinkage priors in high-dimensional linear regression models. Asymptotic results for shrinkage priors have e.g.~been established by \cite{ArmDunLee2013b,SonLia2017}, both papers using the KL and testing strategy. \cite{BaiMorAnt2020,WeiReiHop2020} build upon this work and establish posterior concentration results for splines in sparse additive models. However, all these papers do not consider Bayesian penalized splines with a roughness penalty.

The main practical advantage of Bayesian penalized splines with a roughness penalty is that we do not need to choose the ``right'' number of spline knots to obtain good results in practice. We simply choose a large number of spline knots and let the penalized splines prior with roughness penalty and hyperprior $p(\tau^2)$ on the smoothing variance adaptively prevent overfitting. A different Bayesian strategy to achieve adaptivity is to place a hyperprior on the number of spline knots. This approach has been investigated theoretically by \cite{Jonvan2012}. However, Bayesian penalized splines with a roughness penalty have been shown to be competitive \citep[see, e.g.,][]{LanBre2004} and they are preferable from a computational point of view, especially for non-Gaussian response models. This is because Bayesian penalized splines avoid the need for reversible jump MCMC and they are also compatible with Stan \citep{CarGelHof2017} because of the continuous nature of all involved parameters \citep[see, e.g.,][Section 2.10]{HarRupWan2018}. 

\subsection{Main contributions of this paper}

\begin{itemize}\setlength\itemsep{0em}
    \item We close a gap in the literature and study posterior concentration rates for Bayesian penalized splines. Closing this gap is of great importance as Bayesian penalized splines are widely used in modern statistical analyses. In our main theorem, Theorem~\ref{theo}, we establish sufficient conditions for posterior concentration at near optimal rate.
        \item Based on our theoretical results we can identify specific hyperpriors $p(\tau^2)$ for the smoothing variance that lead to posterior concentration at near optimal rate. These are the first theoretical guarantees in the literature ensuring optimal asymptotic performance for Bayesian penalized splines.
    \item Our theoretical results rely on several new concepts, namely a novel proper version of the partially improper penalized splines prior based on empirical projections as well as an innovative spline estimator. We prove that this {truncated DR estimator} achieves the optimal rate of convergence. 
    \item To define the truncated DR estimator, we establish the concept of DR bases for penalized splines in a rigorous manner. Similar concepts are scattered throughout the literature on penalized splines \cite[e.g.][]{NycCum1996,Rup2002,RupWanCar2003,ClaKriOps2009,Woo2017b}, but a comprehensive introduction has been missing so far.
    \item We complement our theoretical results with empirical evidence demonstrating the adaptivity of the hyperprior on the smoothing variance. Thereby, we use the``scale-dependent'' (Weibull) hyperprior suggested by \cite{KleKne2016}, for which we also establish a new connection to the Generalized Exponential-Type distribution of \cite{ProMab2010}.
\end{itemize}

\subsection{Outline}
This paper is structured as follows: In Section~\ref{sec:PenSplines} we formally introduce Bayesian penalized splines. Section~\ref{sec:tDRB} presents the truncated DR estimator, the key tool for the proof of our main posterior concentration result, Theorem~\ref{theo}. Theorem~\ref{theo} is contained in Section~\ref{sec:PostConc} along with further theoretical results. Section \ref{sec:EmpEvidence} presents the results of a simulation study, where we demonstrate the adaptivity of the hyperprior on the smoothing variance and compare our proper penalized splines prior empirically with the alternative suggestion based on the mixed model reparametrization. 
The final Section~\ref{sec:Discuss} concludes with a discussion. The Supplement contains all our proofs, several auxiliary results and computational details.

\section{Bayesian penalized splines}\label{sec:PenSplines}

Throughout, we consider a univariate Gaussian nonparametric regression model 
\begin{align}\label{GaussiannonparametricRegressionModel}
 Y_i=f_0(x_i)+\varepsilon_i,\ \varepsilon_i\overset{iid}{\sim}N_1(0,\sigma_0^2),\ i=1,\dots,n,
\end{align}
with design points $x_1,\dots,x_n\in[0,1]$, where the aim is to estimate the unknown regression function $f_0:[0,1]\longrightarrow\dsR$. For ease of the exposition, we first assume that the true value $\sigma_0^2\in(0,\infty)$ of the residual variance $\sigma^2$ is known but later we also consider the more realistic case when $\sigma_0^2$ is unknown (see Section~\ref{SecUnknownVar}).

In the frequentist context, \cite{Xia2019} distinguishes between the following three types of penalized splines depending on the particular roughness penalty being used:
\begin{itemize}\itemsep0em
    \item \emph{O'Sullivan penalized splines} or \emph{O-splines} \citep{OS1986,WanOrm2008} use the integrated squared $q$th derivative as roughness penalty.
    \item \emph{P-splines} \citep{EilMar1996} use a particular B-spline basis and 
    the sum of squared finite differences of adjacent B-spline coefficients as roughness penalty.
    \item \emph{T-splines} \citep{Rup2002} use the truncated power basis and a ridge penalty.
\end{itemize}
This distinction is also valid in the Bayesian context, where we have \emph{Bayesian O-splines} \citep{WanOrm2008,HarRupWan2018}, \emph{Bayesian P-splines} \citep{LanBre2004} and \emph{Bayesian T-splines} \citep{BerCarRup2002}. We focus on Bayesian O-splines throughout but our results can easily be carried over to Bayesian P-splines and Bayesian T-splines. 

\subsection{Preliminaries}
Before we introduce the Bayesian O-splines approach in greater detail, we first introduce the vector space of splines and equip it with the required structure. Meanwhile, we introduce our notation for many important quantities required in the remainder of this paper.

To estimate the unknown regression function $f_0$, we use the splines $\mathcal{S}(m,\xi)$, where $m\in\dsN$ is the order and $\xi\in[0,1]^{k_0+2}$ is the vector of spline knots with entries $0=\xi_0<\xi_1<\dots<\xi_{k_0}<\xi_{k_0+1}=1$. Splines are piecewise polynomials that are smoothly joined at the knots, i.e.~for order $m\geq 2$ splines are defined as
\begin{align*}
 \mathcal{S}(m,\xi)=\{f\in C^{m-2}([0,1]):\forall j=0,\dots,k_0: f|_{[\xi_j,\xi_{j+1}]} \in \mathcal{P}_{m-1}\ \},
\end{align*}
where $\mathcal{P}_{m-1}$ is the space of polynomials of degree $m-1$ or less. For order $m=1$ the splines are simply step functions with jumps at the knots. It is well-known that the splines $\mathcal{S}(m,\xi)$ form a vector space of dimension $d=m+k_0$ and a popular basis for the vector space of splines $\mathcal{S}(m,\xi)$ are the B-splines $\{B_1,\dots,B_d\}$.\footnote{To avoid ambiguity, we always refer to the normalized B-splines that arise by repeating the boundary knots $m$ times \cite[see, e.g.,][Section 2.1, for a definition]{Xia2019}.} 

Throughout, we assume that the splines $\mathcal{S}(m,\xi)$ are continuous, i.e.~we assume that $m\geq 2$, so that the splines $\mathcal{S}(m,\xi)$ form a $d$-dimensional subspace of the Banach space $(C([0,1]),\|\cdot\|_\infty)$. We denote the \emph{vector of function evaluations} of a continuous function $f\in C([0,1])$ by $f^n=(f(x_1),\dots,f(x_n))^\top\in\dsR^n$ and the \emph{empirical bilinear form} of two continuous functions $f,g\in C([0,1])$ by $\langle f,g\rangle_n=n^{-1}\sum_{i=1}^n f(x_i)g(x_i)=\langle f^n,g^n\rangle/n,$
where $\langle \cdot,\cdot \rangle$ is the Euclidean scalar product. 

The empirical bilinear form induces the \emph{empirical semi-norm} $\|f\|_n=(\langle f,f \rangle_n)^{1/2}=\|f^n\|_2/\sqrt{n}$, where $\|\cdot\|_2$ is the Euclidean norm. The restriction of the empirical bilinear form to the splines $\mathcal{S}(m,\xi)$ is an inner product if and only if the \emph{Gramian matrix} $G_B=(\langle B_j,B_k\rangle_n)_{j,k=1,\dots,d}$ in terms of the B-splines $\{B_1,\dots,B_d\}$ is positive definite. With the $n\times d$ design matrix $B=(B_1^n,\dots,B_d^n)$ in terms of the B-splines $\{B_1,\dots,B_d\}$, the Gramian matrix is equal to $G_B=B^\top B/n$.

For all $q\in \{1,\dots,m-1\}$ we define the \emph{differential bilinear form of order $q$} on the space of splines $\mathcal{S}(m,\xi)$ by $\langle D^qf,D^q g\rangle_{\mathcal{L}^2}=\int_0^1 D^qf(x)\ D^qg(x) dx$,
where $D^qf$ denotes the $q$th derivative of a spline $f\in \mathcal{S}(m,\xi)$. As the classical derivative of a spline $f\in \mathcal{S}(m,\xi)$ is only well-defined for $q\in\{1,\dots,m-2\}$, we use a weak derivative for the case $q=m-1$. The differential bilinear forms induce the \emph{differential semi-norms} $\lVert D^qf\rVert_{\mathcal{L}^2}=(\langle D^qf,D^qf\rangle_{\mathcal{L}^2})^{1/2}$, $q\in\{1,\dots,m-1\}$. 

We write $R_B^{(q)}=(\langle D^qB_j,D^qB_k \rangle_{\mathcal{L}^2})_{j,k=1,\dots,d}$ for the \emph{roughness penalty matrix} of order $q$ in terms of  B-splines $\{B_1,\dots,B_d\}$ and $Y^n=(Y_1,\dots,Y_n)^\top$ for the vector of observations from  \eqref{GaussiannonparametricRegressionModel}.

\subsection{Bayesian O-splines} 
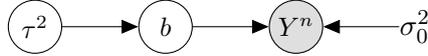
\begin{figure}[H]
  \begin{center}
\begin{tikzpicture}

  \node[obs](y) {$Y^n$};
  \node[latent, left=of y] (b) {${b}$};
  \node[latent, left=of b] (tau2){$\tau^2$};
  \node[const, right=of y] (sigma02) {$\sigma_0^2$};
  
  \edge {tau2} {b} ; 
  \edge {b} {y} ; 
  \edge {sigma02} {y} ; 

\end{tikzpicture}
  \end{center}
  \caption{Graphical model of the O-splines prior in model \eqref{GaussiannonparametricRegressionModel}. Inference about the unknown regression function $f_0$ is based on the marginal posterior $b\mid Y^n,\sigma_0^2$ with the smoothing variance $\tau^2$ integrated out. The marginal posterior $b\mid Y^n,\sigma_0^2$ is generally not available in closed form (see Section~\ref{subsec:ConPostMargPost}).\label{figOSplinesPrior}}
\end{figure}

Next we introduce the Bayesian O-splines approach in greater detail. When using Bayesian O-splines to estimate the unknown regression function $f_0$ in model \eqref{GaussiannonparametricRegressionModel}, we proceed as follows to set up the prior \citep[cf.][]{LanBre2004,WanOrm2008}: 
\begin{enumerate}[a)]
    \item First, we fix the order $m\geq 2$ and the knot vector $\xi$ of the spline space $\mathcal{S}(m,\xi)$. We often use cubic splines of order $m=4$ and about $20-40$ interior spline knots $k_0$ (either equidistant or based on the empirical quantiles of the the design points $x_1,\dots,x_n$).
    \item 
    Next, we fix the B-spline basis $\{B_1,\dots,B_d\}$ of the spline space $\mathcal{S}(m,\xi)$ defined in step a). This allows us to expand any spline $f\in \mathcal{S}(m,\xi)$ uniquely in the form $f=\sum_{j=1}^d B_jb_j$, where $b\in\dsR^d$ is the B-spline coefficient vector. We set up the associated roughness penalty matrix $R_B^{(q)}=(\langle D^qB_j,D^qB_k \rangle_{\mathcal{L}^2})_{j,k=1,\dots,d}$ of order $q\in\{1,\dots,m-1\}$, where $D^qB_j$ is the $q$th (weak) derivative of the B-spline basis function $B_j,\ j=1,\dots,d$. The roughness penalty matrix $R_B^{(q)}$ is a key quantity for Bayesian O-splines as it enables us to quantify the roughness of a spline $f=\sum_{j=1}^d B_jb_j$ in the convenient form
    \begin{align*}
    \|D^qf\|^2_{\mathcal{L}^2}=\int_0^1 (D^q f(x))^2 dx=b^\top R_B^{(q)}b.
    \end{align*}
    A popular choice for the order of the roughness penalty is $q=2$. 
\item 
 Finally, we fix a hyperprior $p(\tau^2)$ on the \emph{smoothing variance} $\tau^2>0$. Then we introduce the \emph{O-splines prior}, which is defined through the following hierarchical specification
 \begin{align}
     b\mid \tau^2&\sim p(b\mid \tau^2) \propto  (\tau^2)^{-(d-q)/2}\ \exp\left(-b^\top R_B^{(q)}b/(2\tau^2)\right), \label{ConditionalOSplinesPrior}\\
     \tau^2&\sim p(\tau^2)\label{Hyperprior}.
 \end{align}
\end{enumerate}
We refer to \eqref{ConditionalOSplinesPrior} as \emph{conditional O-splines prior}. The \emph{marginal O-splines prior} that we effectively place on the B-spline coefficient vector $b\in\dsR^d$ is obtained by integrating the smoothing variance out, i.e.
 \begin{align}\label{OSplinesPrior}
     p(b)\propto \int_0^\infty p(b\mid \tau^2)p(\tau^2)d\tau^2,\ b\in\dsR^d.
 \end{align}
 The O-splines prior prevents overfitting by allocating most of the prior mass to less wiggly splines for which $\|D^qf\|^2_{\mathcal{L}^2}=b^\top R_B^{(q)}b$ is small. Overfitting is an issue due to the relatively large number of spline knots $k_0$ chosen in step a). Figure \ref{figOSplinesPrior} shows a graphical model of the O-splines prior in model~\eqref{GaussiannonparametricRegressionModel}.

The Bayesian O-splines approach with prior \eqref{OSplinesPrior} is the Bayesian counterpart of the frequentist O-splines approach. The \emph{frequentist O-splines estimator} $\widehat{f}_\lambda$ is defined as the minimizer of the penalized sum of squares $\|Y^n-f^n\|^2_2+\lambda\ \lVert D^qf\rVert_{\mathcal{L}^2}^2$ over the space of splines $\mathcal{S}(m,\xi)$, where $\lambda\geq 0$ is a nonnegative \emph{smoothing parameter}. In terms of the B-splines $\{B_1,\dots,B_d\}$, the penalized sum of squares is $\lVert Y^n-Bb\rVert_2^2+\lambda\ b^\top R_B^{(q)}b$ and the frequentist O-splines estimator has the form 
\begin{align}\label{EstimatedCoeffVector}
\widehat{f}_\lambda = \sum_{j=1}^d B_j\widehat{b}_{\lambda,j}\quad\text{with}\quad \widehat{b}_\lambda = \left(B^\top B+\lambda\ R_B^{(q)}\right)^{-1}B^\top Y^n.
\end{align}
The frequentist O-splines estimator~\eqref{EstimatedCoeffVector} has been studied by \cite{OS1986,WanOrm2008,ClaKriOps2009,Xia2019}. 

In the Bayesian O-splines approach we exploit the duality between penalty and prior and replace the roughness penalty $\lVert D^qf\rVert_{\mathcal{L}^2}^2=b^\top R_B^{(q)}b$ by the corresponding prior distribution \eqref{ConditionalOSplinesPrior}. Augmenting the smoothing variance $\tau^2$ with a hyperprior \eqref{Hyperprior}, we arrive at \eqref{OSplinesPrior}. One advantage of the Bayesian approach is that it yields an entire posterior distribution on the splines $\mathcal{S}(m,\xi)$, thus allowing for straightforward uncertainty quantification via posterior credible bands.

\subsection{The hyperprior on the smoothing variance}
In the following, we focus on the hyperprior $p(\tau^2)$ on the smoothing variance $\tau^2$. The hyperprior plays a crucial role for the Bayesian O-splines approach and allows to adapt to the smoothness of the underlying function (see Section~\ref{subsec:Adaptivity} for empirical evidence). 

The distributional family of the hyperprior is a controversial issue and various different families have been suggested in the literature, such as an Inverse Gamma family $\tau^2\sim IG(\alpha,\beta)$, a Gamma family $\tau^2\sim Ga(\alpha,\beta)$, a Weibull family $\tau^2\sim Weibull(k,\lambda)$, a Uniform family $\tau^2\sim U(0,\beta)$ or a Scaled Beta Prime family $\tau^2\sim SBP(\alpha,\beta,\gamma)$ \citep[cf.][]{LanBre2004,KleKne2016,HarRupWan2018}. 

In the following proposition we derive the marginal O-splines prior \eqref{OSplinesPrior} for these suggested hyperpriors on the smoothing variance $\tau^2$. 
\begin{proposition}[Marginal O-splines prior]\label{MarginalPrior} 
\begin{enumerate}[a)]~
      \item Inverse Gamma: Let $\tau^2\sim IG(\alpha,\beta)$, i.e. $p(\tau^2)\propto (\tau^2)^{-\alpha-1}\exp(-\beta/\tau^2),\tau^2 >0$. Then:
      \begin{align*}
    p(b)\propto \left(1+b^\top R_B^{(q)}b/(2\beta)\right)^{-(\alpha+(d-q)/2)},\ b\in\dsR^d.
\end{align*}
      \item Gamma: Let $\tau^2\sim Ga(\alpha,\beta)$, i.e. $p(\tau^2)\propto (\tau^2)^{\alpha-1}\exp(-\beta\tau^2),\tau^2 >0$. Then:
      \begin{align*}
    p(b)\propto \left(\sqrt{b^\top R_B^{(q)}b}\right)^{\alpha-(d-q)/2}K_{\alpha-(d-q)/2}\left(\sqrt{2\beta}\sqrt{b^\top R_B^{(q)}b}\right)\ \dsOne_{\left\lbrace b^\top R_B^{(q)}b>0\right\rbrace},\ b\in\dsR^d,
\end{align*}
where $K(\cdot)$ is the modified Bessel function of the second kind.
\item Weibull: Let $\tau^2 \sim Weibull(k,\lambda)$, i.e. $p(\tau^2)\propto (\tau^2)^{k-1}\ \exp\left(-\lbrace\lambda\tau^2\rbrace^{k}\right),\ \tau^2>0$. Then: 
\begin{align*}
p(b)\propto\mathcal{H}^{2,0}_{0,2}\left( \lambda \dfrac{b^\top R_B^{(q)}b}{2} \Biggm\vert\begin{matrix}  \\ (0,1),(1-(d-q)/(2k),1/k)  \end{matrix}\right)\ \dsOne_{\left\lbrace b^\top R_B^{(q)}b>0\right\rbrace},\ b\in\dsR^d,
\end{align*}
where $\mathcal{H}(\cdot)$ is Fox's H-function.
\item Uniform: Let $\tau^2\sim U(0,\beta)$, i.e. $p(\tau^2)\propto \dsOne_{(0,\beta)}(\tau^2),\tau^2>0.$  Then:
\begin{align*}
    p(b)\propto \left(b^\top R_B^{(q)}b\right)^{-(d-q)/2+1} \Gamma\left((d-q)/2-1,b^\top R_B^{(q)}b/(2\beta)\right)\ \dsOne_{\left\lbrace b^\top R_B^{(q)}b>0\right\rbrace},\ b\in\dsR^d,
\end{align*}
where $\Gamma(\cdot,\cdot)$ is the upper incomplete gamma function.
     \item Scaled Beta Prime \citep[see][for details]{PerPerRam2017}: Let $\tau^2\sim SBP(\alpha,\beta,\gamma)$, i.e. $p(\tau^2)\propto (\tau^2)^{\alpha-1}(\tau^2/\gamma+1)^{-(\alpha+\beta)},\tau^2>0$. Then:
   \begin{align*}
    p(b)\propto U\left(\beta+(d-q)/2,-\alpha +(d-q)/2+1,b^\top R_B^{(q)}b/(2\gamma)\right)\ \dsOne_{\left\lbrace b^\top R_B^{(q)}b>0\right\rbrace},\ b\in\dsR^d,
\end{align*}
where $U(\cdot)$ is the confluent hypergeometric function of the second kind.
\end{enumerate}
\end{proposition}

Proposition~\ref{MarginalPrior} shows that we obtain very different marginal priors \eqref{OSplinesPrior} depending on the hyperprior $p(\tau^2)$ on the smoothing variance. Results a) and b) are familiar results in the penalized splines literature. To the best of our knowledge, Results c)--e) yield new insights that were not recognized so far (see, e.g., \citet{KleKne2016} for Results c), d) and for Result e) see \citet{Sch2011,KleCarKneLanWag2021}, where it is claimed that the respective integrals do not have an analytical solution). We need the indicator in b)--e) to avoid division by zero and to ensure the validity of certain integral representations. However, it does not really impose a restriction as the complementary set $\lbrace b\in\dsR^d: b^\top R_B^{(q)}b=0\rbrace$ is a nullset with respect to $d$-dimensional Lebesgue measure.

\subsection{The conditional and the marginal O-splines posterior}\label{subsec:ConPostMargPost}
Irrespective of the particular choice of the distributional family, the hyperprior $p(\tau^2)$ on the smoothing variance destroys conjugacy and precludes direct access to our posterior of interest. To demonstrate this fact, we next derive the conditional posterior $p(b\mid Y^n,\tau^2,\sigma_0^2)$ and the marginal posterior $p(b\mid Y^n,\sigma_0^2)$. 

The fully Bayesian model specification corresponding to \eqref{GaussiannonparametricRegressionModel}--\eqref{Hyperprior} reads as follows
\begin{align*}
    Y^n\mid b,\sigma_0^2 &\sim N_n(Bb,\sigma_0^2I_n),\\
         b\mid \tau^2&\sim p(b\mid \tau^2) \propto  (\tau^2)^{-(d-q)/2}\ \exp\left(-b^\top R_B^{(q)}b/(2\tau^2)\right),\\
     \tau^2&\sim p(\tau^2).
\end{align*}
It is straightforward to show that the conditional posterior $b\mid Y^n,\tau^2,\sigma_0^2$ of the B-spline coefficients $b\in\dsR^d$ is a Gaussian distribution $N_d(\mu,\Sigma)$ with moments 
\begin{align}\label{ConditionalPosteriorMoments}
\mu=\Sigma/\sigma_0^2\ B^\top Y^n    \quad \text{and}\quad  \Sigma=\left(B^\top B/\sigma_0^2+R_B^{(q)}/\tau^2\right)^{-1}.
\end{align}

Next, we derive the marginal posterior. The marginal posterior of the B-spline coefficients $b\in\dsR^d$ is proportional to
\begin{align}\label{MarginalPosterior}
p(b\mid Y^n,\sigma_0^2)&\propto p(Y^n\mid b,\sigma_0^2)\ p(b)  \nonumber \\
 &\propto N_n(Y^n;Bb,\sigma_0^2I_n)\ p(b)\nonumber\\
 &\propto \exp\left(-\dfrac{1}{2}\left\lbrace b^\top \dfrac{B^\top B}{\sigma_0^2}b-2b^\top \dfrac{B^\top Y^n}{\sigma_0^2}\right\rbrace \right)\ p(b),\ b\in\dsR^d,
\end{align}
where $N_n(Y^n;Bb,\sigma_0^2I_n)$ in the second line denotes the density of a $n$-variate Gaussian distribution with mean vector $Bb$ and covariance matrix $\sigma_0^2I_n$ evaluated at the vector of observations $Y^n\in\dsR^n$. For all hyperpriors $p(\tau^2)$ considered in Proposition~\ref{MarginalPrior}, the marginal posterior \eqref{MarginalPosterior} does not correspond to a familiar distribution. 

In summary it holds: While the conditional posterior $b\mid Y^n,\tau^2,\sigma_0^2$ in model \eqref{GaussiannonparametricRegressionModel}
is Gaussian with known moments \eqref{ConditionalPosteriorMoments}, the marginal posterior $b\mid Y^n,\sigma_0^2$ is generally intractable. This is problematic as inference about the unknown regression function $f_0$ in model \eqref{GaussiannonparametricRegressionModel} is based on the marginal posterior $b\mid Y^n,\sigma_0^2$ and not on the conditional posterior $b \mid Y^n,\tau^2,\sigma_0^2$ \citep[cf.][]{FahKneKon2010}.

From a practical point of view, the loss of conjugacy is not a big issue as we can use Markov chain Monte Carlo (MCMC) methods to sample from the marginal posterior \eqref{MarginalPosterior}. To this end we can e.g.~use a Gibbs sampler, where we draw in an alternating manner from the full conditional posteriors $b\mid Y^n,\tau^2,\sigma_0^2$ and $\tau^2\mid Y^n,b,\sigma_0^2$ \citep[see, e.g.,][]{LanBre2004}. Thereby, we obtain a sample from the joint posterior $(b,\tau^2)\mid Y^n,\sigma_0^2$ and the embedded sample of B-spline coefficients $b^{(1)},\dots,b^{(L)}$ follows the marginal posterior \eqref{MarginalPosterior}. In standard fashion, the sampled coefficients allow us to estimate the unknown regression function $f_0$ via the marginal posterior mean $\widehat{f} = \sum_{j=1}^d B_j\widehat{b_j}$ with $\widehat{b} =\dsE\left[b\mid Y^n,\sigma_0^2\right] \approx 1/L \sum_{l=1}^L b^{(l)}$,
and to quantify uncertainty via posterior credible bands. 

However, the loss of conjugacy caused by the hyperprior $p(\tau^2)$ on the smoothing variance and the resulting intractability of the marginal posterior \eqref{MarginalPosterior} pose great challenges for the asymptotic analysis. This is because we cannot use a ``direct approach'', where we show that the posterior mean converges to the truth and that the posterior variance vanishes sufficiently fast. This renders the asymptotic analysis much more complicated \citep[cf.][Chapter 8]{Ghovan2017}.

A general strategy to establish posterior concentration rates that does not require analytic expressions for the posterior moments is the KL and testing strategy \citep[see, e.g.,][Section 3]{Rou2016}. However, the KL and testing strategy requires a proper prior and the O-splines prior is improper, i.e.~for both the conditional O-splines prior \eqref{ConditionalOSplinesPrior} and the marginal O-splines prior \eqref{OSplinesPrior} it holds
\begin{align}\label{Impropriety}
    \int_{\dsR^d} p(b\mid \tau^2) db=+\infty\quad \text{and}\quad \int_{\dsR^d} p(b) db=+\infty.
\end{align}
The impropriety~\eqref{Impropriety} is due to the fact that the roughness penalty matrix $R_B^{(q)}$ is rank-deficient, which in turn results from the fact that the $q$th derivative $D^q$ vanishes on the space of polynomials $\mathcal{P}_{q-1}$. We speak of \emph{partial impropriety} because both priors \eqref{ConditionalOSplinesPrior} and \eqref{OSplinesPrior} decompose into the product of an improper prior and a proper prior (see explanation below Proposition~\ref{propPropertiesDRB}). 
While the partial impropriety is typically not an issue for posterior sampling and the propriety of the posterior \citep[see, e.g.,][]{HenBreFah2006,FahKne2009,KleKne2016}, a proper prior is necessary to be able to use the KL and testing strategy. 

\subsection{The proper O-splines prior}
In this section we introduce a proper version of the partially improper O-splines prior. This \emph{proper O-splines prior} allows us to use the KL and testing strategy, which enables us to establish posterior concentration rates for Bayesian penalized splines despite the analytical intractability of the marginal posterior \eqref{MarginalPosterior}. 

To turn the usual O-splines prior into a proper prior, we introduce an additional penalty for the polynomials in $\mathcal{P}_{q-1}$ and make the following definition.

\begin{definition}[Proper O-splines prior]
Let $X_q$ denote the full rank $n\times q$ design matrix of the monomial basis $\{1,x,\dots,x^{q-1}\}$ evaluated at the design points $x_1,\dots,x_n\in [0,1]$. Let $H^{(q)}=X_q(X_q^\top X_q)^{-1}X_q^\top$ denote the projection matrix from $\dsR^n$ onto $\text{span}(X_q)$. Then we define:

\begin{enumerate}[a)]
    \item The \emph{proper conditional O-splines prior of order $q\in\{1,\dots,m-1\}$} on the B-spline coefficients $b\in \dsR^d$ is defined as the multivariate Gaussian density
\begin{align}\label{ConditionalProperOSplinesPrior}
   \widetilde{p}(b\mid \tau^2)=N_d\left(b;0,\left({R_B^{(q)}}/{\tau^2}+{B^\top H^{(q)}B}/{(n\tau^2_{poly})}\right)^{-1}\right),\ b\in\dsR^d,
\end{align}
where $\tau^2>0$ is the smoothing variance as before and $\tau^2_{poly}>0$ is an additional hyperparameter that controls the strength of the penalty for the polynomial part.  

\item The \emph{proper marginal O-splines prior of order} $q\in\{1,\dots,m-1\}$ on the B-spline coefficients $b\in \dsR^d$ is defined as the density
\begin{align}\label{ProperOSplinesPrior}
    \widetilde{p}(b)=\int_0^\infty \widetilde{p}(b\mid \tau^2)p(\tau^2)d\tau^2,\ b\in\dsR^d.
\end{align}
\end{enumerate}
\end{definition}
We regard \eqref{ConditionalProperOSplinesPrior} as a proper approximation of the conditional O-splines prior \eqref{ConditionalOSplinesPrior} and \eqref{ProperOSplinesPrior} as a proper approximation of the marginal O-splines prior \eqref{OSplinesPrior}. In Section~\ref{subsec:ProperPriors} we will show that this view is justified empirically.  

Next, we shift our perspective from the prior and posterior on the B-spline coefficients $b\in\dsR^d$ to the corresponding prior and posterior on the splines $\mathcal{S}(m,\xi)$, i.e.~we shift our perspective from probability measures on Euclidean space $\dsR^d$ to probability measures on function spaces. The proper O-splines prior \eqref{ProperOSplinesPrior} {induces} a prior distribution $\widetilde{\Pi}$ on the splines $\mathcal{S}(m,\xi)$. Loosely speaking, the induced prior $\widetilde{\Pi}$ is simply the law of the random spline function $f=\sum_{j=1}^d B_jb_j$ when $b\in\dsR^d$ follows density \eqref{ProperOSplinesPrior}. Formally, $\widetilde{\Pi}$ is a probability measure on the Borel sets $\mathcal{B}_C$ of the continuous functions $C([0,1])$ with support $\mathcal{S}(m,\xi)$ that arises as {pushforward} under the Borel-measurable B-spline expansion $b\mapsto \sum_{j=1}^d B_jb_j$, i.e.
\begin{align}\label{ProperInducedPrior}
    \widetilde{\Pi}(A)=\int_{\{b\in\dsR^d:f=\sum_{j=1}^d B_jb_j\in A\}} \widetilde{p}(b)db,\ A\in\mathcal{B}_C.
\end{align}

Given observations $Y^n\in\dsR^n$ from model \eqref{GaussiannonparametricRegressionModel}, the (marginal) O-splines posterior $\widetilde{\Pi}(\cdot\mid Y^n,\sigma_0^2)$ on the splines $\mathcal{S}(m,\xi)$, the key quantity of interest in the context of posterior concentration rates, is then obtained by the following nonparametric version of Bayes' rule
 \begin{align}\label{BayesRuleFunctions}
     \widetilde{\Pi}(A\mid Y^n,\sigma_0^2)=\dfrac{\int_A N_n(Y^n;f^n,\sigma_0^2I_n)\ \widetilde{\Pi}(df)}{\int_{C([0,1])}N_n(Y^n;f^n,\sigma_0^2I_n) \ \widetilde{\Pi}(df)},\ A\in \mathcal{B}_C.
 \end{align} 
 For further technical details on priors, posteriors and Bayes' rule in the nonparametric setting we refer to \citet{Ghovan2017}. 
 
 In Section~\ref{sec:PostConc} we construct a sequence of priors of the form~\eqref{ProperInducedPrior} that leads to a sequence of posteriors of the form~\eqref{BayesRuleFunctions} by Bayes' rule. In our main posterior concentration result, Theorem~\ref{theo}, we establish sufficient conditions which guarantee that this sequence of posteriors concentrates around the unknown regression function $f_0$ at near optimal rate. In the following Section~\ref{sec:tDRB} we introduce our key tool for the proof of Theorem~\ref{theo},  the \textit{truncated DR estimator}.

\section{The truncated DR estimator}\label{sec:tDRB}

In this section we introduce the truncated DR estimator $\widehat{f}_t$, the key tool for the proof of our main posterior concentration result, Theorem~\ref{theo}, which is stated in the subsequent Section~\ref{sec:PostConc}. 
The truncated DR estimator $\widehat{f}_t$ is closely related to the frequentist O-splines estimator $\widehat{f}_\lambda$ \eqref{EstimatedCoeffVector} and is motivated as follows:
\begin{itemize}
    \item The main difficulty when using the KL and testing strategy in the present setting is to devise a suitable test $\phi$ for the test problem $H_0:f=f_0$ vs. $H_1:\|f-f_0\|_n\geq M\epsilon_n$ \cite[cf.][Section 3]{Rou2016}. A test $\phi=\phi(Y^n)$ is suitable if its type I error probability vanishes and if its type II error probability vanishes sufficiently fast.
    \item A natural candidate for this test is $\phi=\dsOne\{\|\widehat{f}_\lambda-f_0\|_n\geq M/2\ \epsilon_n\}$, where $\widehat{f}_\lambda$ is the frequentist O-splines estimator \eqref{EstimatedCoeffVector}. Using Markov's inequality and the results of \cite{ClaKriOps2009,Xia2019}, one can show that the type I error probability of this test vanishes. However, the type II error probability of the test based on $\widehat{f}_\lambda$ is very difficult to control. 
    \item Therefore, we introduce the truncated DR estimator $\widehat{f}_t$. We will show that the type I error probability of the corresponding test $\phi=\dsOne\{\|\widehat{f}_t-f_0\|_n\geq M/2\ \epsilon_n\}$ vanishes and that the type II error probability of the test based on $\widehat{f}_t$ is much easier to control than for the test based on the frequentist O-splines estimator $\widehat{f}_\lambda$.
\end{itemize}

\subsection{The DR bases for O-splines}
The truncated DR estimator $\widehat{f}_t$ projects the observations $Y^n$ from model \eqref{GaussiannonparametricRegressionModel} onto the first basis functions of a {DR basis}. Thus, we first need to introduce the DR bases before we can introduce the estimator itself. 
 
The DR bases are special orthonormal bases of the space of splines $\mathcal{S}(m,\xi)$. These bases were introduced by \cite{DemRei1975} in the context of smoothing splines \cite[see also][Section 2]{Spe1985}. In the following, we carry their idea over from smoothing splines to O-splines. 

O-splines and smoothing splines are closely related but in contrast to smoothing splines the knots $\xi$ and the design points $x_1,\dots,x_n$ as well as the order of the spline $m$ and the order of the roughness penalty $q$ are decoupled for O-splines, which provides additional flexibility and considerable computational savings \citep{WanOrm2008,Woo2017}. We are not the first to carry the idea of \cite{DemRei1975} over from smoothing splines to penalized splines \cite[e.g.][]{NycCum1996,Rup2002,RupWanCar2003,ClaKriOps2009,Woo2017b}. However, a rigorous introduction of the notion of a DR basis for O-splines seems to be missing. We make the following definition.

\begin{definition}[DR basis]\label{DefinitionDRB} Let $\{Z_1,\dots,Z_d\}$ be a basis of the space of splines $\mathcal{S}(m,\xi)$ of order $m\geq 2$. We say that $\{Z_1,\dots,Z_d\}$ is a \emph{DR basis of order} $q\in\{1,\dots,m-1\}$ if it is orthonormal with respect to the empirical bilinear form and orthogonal with respect to the differential bilinear form of the corresponding order, i.e. $\{Z_1,\dots,Z_d\}$ is a \emph{DR basis of order} $q\in\{1,\dots,m-1\}$ if 
    \begin{align*}
     \langle Z_j,Z_k\rangle_n=\delta_{j,k}\quad\text{ and }\quad\langle D^qZ_j,D^qZ_k\rangle_{\mathcal{L}^2}=\gamma_j^{(q)}\delta_{j,k},\ j,k=1,\dots,d,
    \end{align*}
where $\delta_{j,k}$ denotes Kronecker's delta and $0\leq \gamma_1^{(q)}\leq \dots \leq \gamma_d^{(q)}$ are nondecreasing nonnegative numbers, which we refer to as \emph{DR eigenvalues}.
\end{definition}

\begin{proposition}[Existence]\label{propExistenceDRB} Let $\{B_1,\dots,B_d\}$ be the B-spline basis of the spline space $\mathcal{S}(m,\xi)$ of order $m\geq 2$ and $x_1,\dots,x_n$ some design points in $[0,1]$. Then it holds: There exists a DR basis $\{Z_1,\dots,Z_d\}$ of order ${q\in\{1,\dots,m-1\}}$ if and only if the Gramian matrix $G_B=(\langle B_j,B_k\rangle_n)_{j,k=1,\dots,d}=B^\top B/n$ in terms of the B-splines is positive definite.
\end{proposition}

\begin{figure}[htbp]
\begin{center}
    \includegraphics[width=\textwidth]{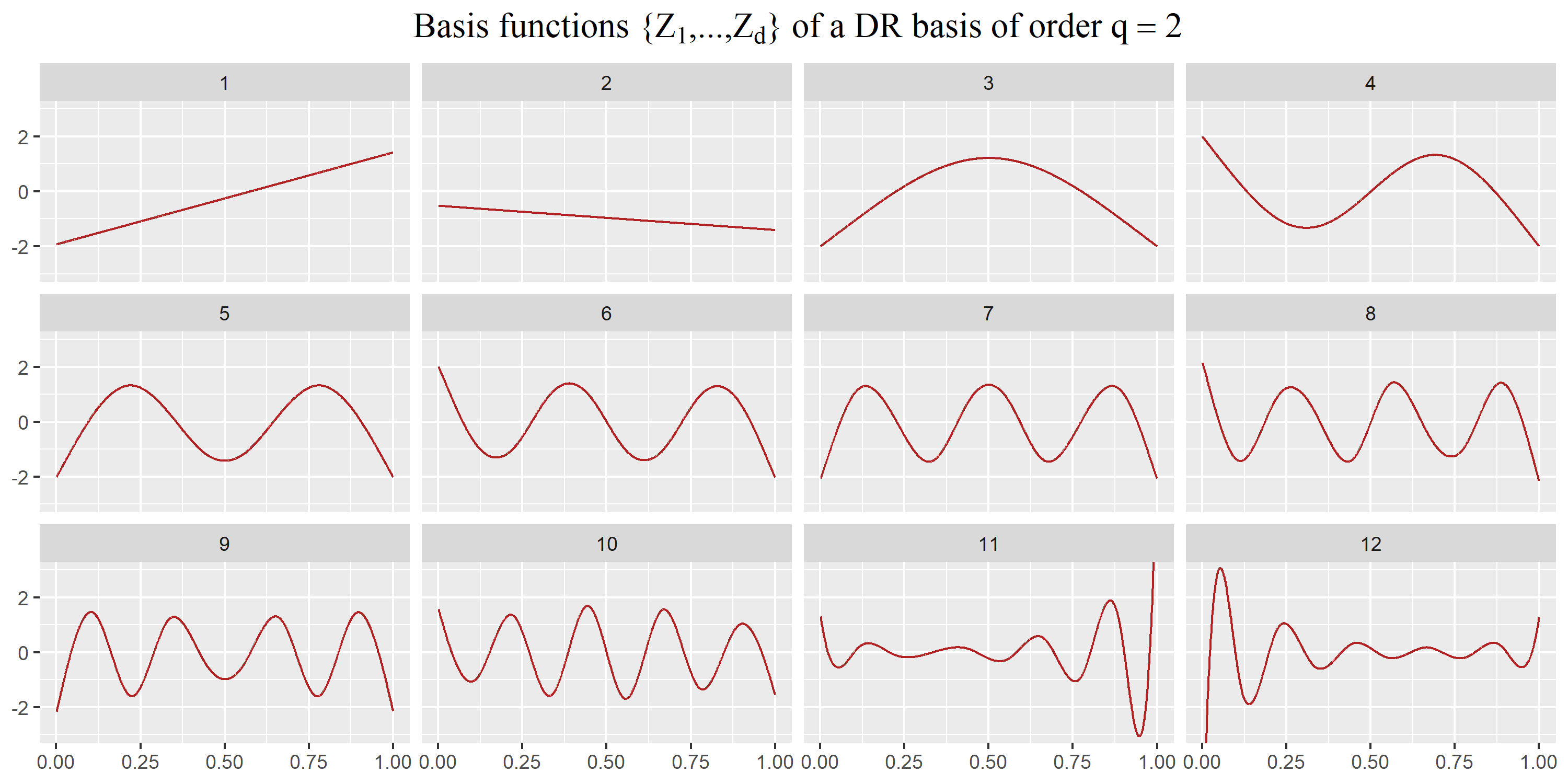}
    \end{center}
\caption{Shown are the basis functions $\{Z_1,\dots,Z_d\}$ of a DR basis of order $q=2$ for the $d=12$ dimensional space of cubic splines $\mathcal{S}(4,\xi)$ with $k_0=8$ equidistant interior knots in the interval $[0,1]$. The $n=1000$ design points $x_1,\dots,x_n$ were chosen in a regular manner, i.e. $x_i=i/n,\ i=1,\dots,n$. The R package \texttt{geigen} \citep{Has2019} was used to compute the transition matrix $A^{(q)}$ in \eqref{ConstructionDRB} by solving the generalized symmetric eigenvalue problem for the pair $(R_B^{(q)},G_B)$. The code for $R_B^{(q)}$ was taken from the Appendix of \cite{WanOrm2008}. By definition~\ref{DefinitionDRB}, the basis functions $\{Z_1,\dots,Z_d\}$ with design matrix $Z=BA^{(q)}$ are orthonormal, i.e.~the Gramian matrix $G_Z=(\langle Z_j,Z_k\rangle_n)_{j,k=1,\dots,d}=Z^\top Z/n$ is equal to the identity $I_d$. Moreover, the roughness penalty matrix $R_Z^{(q)}=(\langle D^qZ_j,D^qZ_k \rangle_{\mathcal{L}^2})_{j,k=1,\dots,d}$ is diagonal with DR eigenvalues $0\leq \gamma_1^{(q)}\leq \dots \leq \gamma_d^{(q)}$ as diagonal entries.}\label{fig2}
\end{figure}

Proposition~\ref{propExistenceDRB} proves that a DR basis of order $q\in\{1,\dots,m-1\}$ exists if (and only if) the Gramian matrix $G_B=B^\top B/n$ in terms of the B-splines $\{B_1,\dots,B_d\}$ is positive definite or, equivalently, if the B-spline design matrix $B$ of size $n\times d$ has rank $d$.  In this case, we can construct a DR basis by solving the generalized symmetric eigenvalue problem \citep[see, e.g.,][Chapter 15]{Par1998} for the ordered pair of matrices $(R_B^{(q)},G_B)$. Denoting the resulting $d\times d$ matrix of generalized eigenvectors by $A^{(q)}$, the basis $\{Z_1,\dots,Z_d\}$ with
\begin{align}\label{ConstructionDRB}
    Z_j=\sum_{k=1}^d B_k A_{kj}^{(q)},\ j=1,\dots,d,
\end{align}
is a DR basis of order $q$. The values $\gamma_{1}^{(q)},\dots,\gamma_{d}^{(q)}$ are the generalized  eigenvalues of the pair $(R_B^{(q)},G_B)$, which explains why we refer to them as DR eigenvalues. 
Figure~\ref{fig2} shows the basis functions $\{Z_1,\dots,Z_d\}$ of a DR basis of order $q=2$.

The truncated DR estimator $\widehat{f}_t$ divides the DR basis functions $\{Z_1,\dots,Z_d\}$ shown in Figure~\ref{fig2} into low-frequency functions $\{Z_1,\dots,Z_t\}$ and high-frequency functions $\{Z_{t+1},\dots,Z_d\}$ and then projects the observations $Y^n$ onto the low-frequency functions only. Before we formally introduce the truncated DR estimator, we first address the important question of uniqueness and summarize convenient properties of DR bases for O-splines.

\begin{proposition}[Uniqueness]\label{propUniquenessDRB}
 Let $\mathcal{S}(m,\xi)$ be a spline space of order $m\geq 2$ and let $\{Z_1,\dots,Z_d\}$ be a DR basis of order $q\in\{1,\dots,m-1\}$ with DR design matrix $Z$ of size $n\times d$ and DR eigenvalues $\gamma_{1}^{(q)},\dots,\gamma_{d}^{(q)}$. Then it holds: 
\begin{enumerate}[i)]
 \item The DR eigenvalues are uniquely specified: If $\{\widetilde{Z}_1,\dots,\widetilde{Z}_d\}$ is another DR basis of the same order $q$ with DR eigenvalues $\widetilde{\gamma}_{1}^{(q)},\dots,\widetilde{\gamma}_{d}^{(q)}$, then $\widetilde{\gamma}_j^{(q)}=\gamma_j^{(q)}$ for all $j=1,\dots,d$. Furthermore, the transition matrix from $\{Z_1,\dots,Z_d\}$ to $\{\widetilde{Z}_1,\dots,\widetilde{Z}_d\}$ is blockdiagonal with orthogonal blocks whose size is equal to the multiplicities of the DR eigenvalues.
 \item The first $q$ DR eigenvalues are equal to $0$, while the remaining $d-q$ DR eigenvalues are positive, i.e.~ $0=\gamma_1^{(q)}=\dots=\gamma_q^{(q)}<\gamma_{q+1}^{(q)}\leq \dots \leq \gamma_d^{(q)}$. Moreover, the first $q$ basis functions $\{Z_{1},\dots,Z_q\}$ form an orthonormal basis of the polynomials $\mathcal{P}_{q-1}$, while the remaining $d-q$ basis functions $\{Z_{q+1},\dots,Z_d\}$ form an orthonormal basis of the orthogonal complement of the space $\mathcal{P}_{q-1}$ in the space of splines $\mathcal{S}(m,\xi)$. 
 \end{enumerate}
 \end{proposition}
 
Proposition~\ref{propUniquenessDRB} shows that the DR eigenvalues $\gamma_{1}^{(q)},\dots,\gamma_{d}^{(q)}$ are uniquely specified. In contrast to that, the basis functions of a DR basis of order $q\in\{1,\dots,m-1\}$ in the sense of Definition~\ref{DefinitionDRB} are not uniquely specified. 
The first $q$ basis functions are an arbitrary orthonormal basis of the space of polynomials $\mathcal{P}_{q-1}$. Moreover, the remaining $d-q$ basis functions are only unique (up to sign) if the positive DR eigenvalues are distinct, i.e.~if $\gamma_{q+1}^{(q)}<\dots<\gamma_d^{(q)}$.   

However, the non-uniqueness of the DR basis functions is not an issue for us as the precise choice of the DR basis is irrelevant for our purpose to establish posterior concentration rates. We can simply choose any DR basis of order $q\in\{1,\dots,m-1\}$ in the sense of Definition \ref{DefinitionDRB}. All these bases of $\mathcal{S}(m,\xi)$ share the same relevant properties, which are summarized in the following proposition.

\begin{proposition}[Properties of a DR basis]\label{propPropertiesDRB} Let $\mathcal{S}(m,\xi)$ be a spline space of order $m\geq 2$ and let $\{Z_1,\dots,Z_d\}$ be a DR basis of order $q\in\{1,\dots,m-1\}$ with DR design matrix $Z$ of size $n\times d$ and DR eigenvalues $\gamma_{1}^{(q)},\dots,\gamma_{d}^{(q)}$. Then it holds:
\begin{enumerate}[i)]
\item We can express any spline $f\in\mathcal{S}(m,\xi)$ uniquely in the form $f=\sum_{j=1}^d Z_ju_j$ for some coefficient vector $u\in\dsR^d$.
The squared empirical norm and the integrated squared $q$th derivative of $f=\sum_{j=1}^d Z_ju_j$ are $\lVert f\rVert_n^2=\sum_{j=1}^d u_j^2$ and $\lVert D^qf\rVert_{\mathcal{L}^2}^2=\sum_{j=q+1}^d \gamma_j^{(q)}u_j^2$, respectively.
 \item In terms of the DR basis $\{Z_1,\dots,Z_d\}$, the frequentist O-splines estimator \eqref{EstimatedCoeffVector} has the form
 \begin{align}\label{OSplinesEstimatorDRB}
\widehat{f}_\lambda=\sum_{j=1}^d Z_j w_{\lambda,j}\widehat{u}_{j},     
 \end{align}
where $\widehat{u}=(Z^\top Z)^{-1}Z^\top Y^n=Z^\top Y^n/n$ is the least squares coefficient vector and the $w_{\lambda,j},\ j=1,\dots,d,$ are nonincreasing ``shrinkage weights'' of the form $w_{\lambda,j}=\left(1+\lambda/n\ \gamma_j^{(q)}\right)^{-1}$ for $j=1,\dots,d,$ that lie in the interval $(0,1].$
 \item The proper conditional O-splines prior \eqref{ConditionalProperOSplinesPrior} in terms of the coefficient vector $u\in\dsR^d$ is
  \begin{align}\label{ProperOSplinesPriorDemmlerReinsch}
     \widetilde{p}(u\mid \tau^2) = \prod_{j=1}^q N_1\left(u_j;0,\tau^2_{\mbox{\scriptsize{poly}}}\right)\ \prod_{j=q+1}^d N_{1}\left(u_j;0,\tau^2 /\gamma_{j}^{(q)}\right).
 \end{align}
\end{enumerate} 
\end{proposition}

Proposition \ref{propPropertiesDRB} shows that the DR bases are very useful in the context of O-splines because they simplify the mathematical expressions of important quantities and provide further insight. Formula \eqref{ProperOSplinesPriorDemmlerReinsch}, for instance, shows that the proper conditional O-splines prior \eqref{ConditionalProperOSplinesPrior} has a very simple form in terms of the coefficients $u\in\dsR^d$ of a DR basis $\{Z_1,\dots,Z_d\}$. Moreover, \eqref{ProperOSplinesPriorDemmlerReinsch} reveals the effect of the additional penalty ${B^\top H^{(q)}B}/{(n\tau^2_{poly})}$ that we have introduced for our definition of the proper conditional O-splines prior \eqref{ConditionalProperOSplinesPrior} -- without this additional penalty, the prior on the first $q$ polynomial coefficients $u_j,\ j=1,\dots,q,$ is improper (explaining the notion of partial impropriety introduced earlier).

\subsection{The truncated DR estimator}
 
We are now able to define the truncated DR estimator $\widehat{f}_t$, which is inspired by the form \eqref{OSplinesEstimatorDRB} of the frequentist O-splines estimator $\widehat{f}_\lambda$ in terms of a DR basis.
\begin{definition}[Truncated DR estimator]\label{DefinitiontDRB} Let $\mathcal{S}(m,\xi)$ be a spline space of order $m\geq 2$ and let $\{Z_1,\dots,Z_d\}$ be a DR basis of order $q\in\{1,\dots,m-1\}$ with DR design matrix $Z$ of size $n\times d$. Let further $\widehat{u}=Z^\top Y^n/n$ the corresponding least squares coefficient vector and let $t\in \{1,\dots,d\}$ a \emph{cut-off} or \emph{truncation point}. Let $w_t=(1,\dots,1,0,\dots,0)^\top$ with $w_{t,j}=1$ for $j=1,\dots,t,$ and $w_{t,j}=0$ for $j=t+1,\dots,d$. Then we call an estimator of the form
\begin{align}\label{DefTruncatedDR}
    \widehat{f}_t=\sum_{j=1}^d Z_j w_{t,j}\widehat{u}_{j}=\sum_{j=1}^t Z_{j}\widehat{u}_j
\end{align}
a \emph{truncated DR estimator}.
\end{definition}

The following proposition reveals the main advantage of a truncated DR estimator $\widehat{f}_t$ over the O-splines estimator $\widehat{f}_\lambda$ in the present context: To control the type II error probability of the test $\phi=\dsOne\{\|\widehat{f}-f_0\|_n\geq M/2\ \epsilon_n\}$, we essentially need to control the distribution of the squared empirical distance $\|\widehat{f}-f\|_n^2$ when the observations $Y^n$ are distributed as $Y^n\sim N_n(f^n,\sigma_0^2I_n)$ where $f\in \mathcal{S}(m,\xi)$ is a spline function. 

\begin{proposition}[Type II error]\label{propType2err} Let $\mathcal{S}(m,\xi)$ be a spline space of order $m\geq 2$ and let $\{Z_1,\dots,Z_d\}$ be a DR basis of order $q\in\{1,\dots,m-1\}$ with DR design matrix $Z$. Let $\widehat{f}_t=\sum_{j=1}^t Z_j\widehat{u}_j$ be the corresponding truncated DR estimator with cut-off $t\in\{1,\dots,d\}$ and $\widehat{u}=Z^\top Y^n/n$. Let further $\widehat{f}_\lambda$ be the O-splines estimator and let $f=\sum_{j=1}^d Z_ju_j$ be a spline function in $\mathcal{S}(m,\xi)$. Then it holds under $Y^n\sim N_n(f^n,\sigma_0^2I_n):$ 
\begin{align}\label{DistriAlternative}
    \|\widehat{f}_\lambda-f\|_n^2\sim {\sigma_0^2}/{n}\ \sum_{j=1}^d w_{\lambda,j}^2  \chi^2_j\quad\text{and}\quad   \|\widehat{f}_t-f\|_n^2\sim {\sigma_0^2}/{n}\  \chi^2(t)+\sum_{j=t+1}^d u_j^2,
\end{align}
where the $\chi_1^2,\dots,\chi_d^2$ are independent noncentral chi-square distributed variates with $1$ degree of freedom and noncentrality parameters $nu_j^2(w_{\lambda,j}-1)^2/(\sigma_0^2w_{\lambda,j}^2)$ and the $\chi^2(t)$ on the right is a single (central) chi-square distributed variate with $t$ degrees of freedom.
\end{proposition}

As the distribution of $\|\widehat{f}_t-f\|_n^2$ on the right of \eqref{DistriAlternative} is much easier to handle than the distribution of $\|\widehat{f}_\lambda-f\|_n^2$ on the left of \eqref{DistriAlternative}, the type II error probability for the test based on $\widehat{f}_t$ is much easier to control than for the test based on $\widehat{f}_\lambda$. 

\subsection{Remarks}
\begin{itemize}
    \item A truncated DR estimator $\widehat{f}_t$ divides the DR basis functions $\{Z_1,\dots,Z_d\}$ into low-frequency basis functions $\{Z_1,\dots,Z_t\}$ and high-frequency basis functions $\{Z_{t+1},\dots,Z_d\}$ and then projects the observations $Y^n$ onto the low-frequency basis functions only. 
    \item The relation between a truncated DR estimator $\widehat{f}_t$ and the O-splines estimator $\widehat{f}_\lambda$ is comparable to the relation of principal components regression and the usual ridge regression estimator \citep[see, e.g., Section 3.5.1 of][]{HasTibFri2009}.
    \item Conceptually, a truncated DR estimator is also comparable to a truncated Fourier series estimator \cite[][]{EubHarSpe1990} as we project the observations onto the first basis functions of an orthonormal basis that is ordered in terms of increasing complexity. However, there are large differences in the details. 
    \item In Proposition~C.2 in the Supplement we show that a truncated DR estimator $\widehat{f}_t$ is able to achieve the optimal rate of convergence in model \eqref{GaussiannonparametricRegressionModel}. This is an important intermediate result for the proof of our main posterior concentration result, Theorem~\ref{theo}, as it implies that the type I error probability of the test $\phi=\dsOne\{\|\widehat{f}_t-f\|_n\geq M/2\ \epsilon_n\}$ vanishes. 
\end{itemize} 

In the following Section~\ref{sec:PostConc} we investigate the asymptotic behavior of the Bayesian O-splines approach. Thereby, we use the test $\phi=\dsOne\{\|\widehat{f}_t-f_0\|_n\geq M/2\ \epsilon_n\}$ based on a truncated DR estimator $\widehat{f}_t$ to establish posterior concentration rates.

\newpage

\section{Posterior concentration rates for Bayesian O-splines}\label{sec:PostConc}
In this section we state our main posterior concentration result for Bayesian penalized splines, Theorem~\ref{theo}. Before we state the theorem, we first define a sequence of priors, which leads to a sequence of posteriors by Bayes' rule, and state our assumptions.

To investigate the asymptotic behavior of the O-splines approach with prior \eqref{ProperOSplinesPrior} on the B-spline coefficients $b\in \dsR^d$, we construct a sequence $\widetilde{\Pi}_n,\ n\in\mathbb{N},$ of proper induced O-splines priors \eqref{ProperInducedPrior}, which is supported on a sequence of spline spaces $\mathcal{S}(m,\xi_n),\ n\in \mathbb{N},$ of increasing dimension $d=d(n)$. Given observations $Y^n\in\dsR^n$ from model \eqref{GaussiannonparametricRegressionModel}, we obtain a sequence of O-splines posteriors $\widetilde{\Pi}_n(\cdot\mid Y^n,\sigma_0^2),\ n\in\mathbb{N},$ via Bayes' rule
\begin{align*}
     {\widetilde{\Pi}}_n(A\mid Y^n,\sigma_0^2)=\dfrac{\int_A N_n(Y^n;f^n,\sigma_0^2I_n)\ {\widetilde{\Pi}_n}(df)}{\int_{C([0,1])}N_n(Y^n;f^n,\sigma_0^2I_n) \ {\widetilde{\Pi}_n}(df)},\ A\in \mathcal{B}_C.
 \end{align*}
 
 Ideally, this sequence of posteriors concentrates around the true regression function $f_0$ in the sense that shrinking neighbourhoods of $f_0$ receive more and more posterior mass as the sample size $n$ grows. This idea is formalized by the notion of posterior concentration rates \citep[see, e.g.,][]{GhoGhovan2000,Ghovan2007,Rou2016}. Whether the O-splines posterior concentrates at all around $f_0$ as well as the rate of concentration $\epsilon_n$ if it does, depend crucially on our assumptions, which we state next. 

 \subsection{Assumptions}
 We make the following assumptions for the true model \eqref{GaussiannonparametricRegressionModel}, the parameters of our prior $\widetilde{\Pi}_n=\widetilde{\Pi}_n(m,\xi_n,q,p_n(\tau^2),\tau^2_{poly})$ and the design points $x_1,\dots,x_n\in[0,1]$:
  \begin{itemize}
   \setlength\itemsep{1em}
 \item[A1)] The true regression function $f_0$ is in $C^{m_0}([0,1])$ for some $m_0\in\dsN$ and the residual variance $\sigma_0^2$ is equal to some positive constant (treated as known).
 \item[A2)] The order $m$ of the spline spaces $\mathcal{S}(m,\xi_n)$ is a constant integer and exceeds the regularity of the unknown regression function $f_0$, i.e. $m>m_0$. 
   \item[A3)] The number of interior spline knots $k_0$ diverges at rate $k_0=k_0(n)\sim c_{k_0}n^\delta$, where $c_{k_0}>0$ is some positive constant and $1/(2m_0+1)<\delta<1$. The spline knots $\xi_n$ are equidistant, i.e. $\xi_{nj}=jh$ for $j=0,\dots,k_0+1,$ with mesh width $h=h(n)=1/(k_0+1)$. 
   \item[A4)] The order of the derivative $q$ in the roughness penalty is a constant integer and matches the regularity of the unknown regression function $f_0$, i.e. $q=m_0$.
  \item[A5)] The smoothing variance $\tau^2$ is equipped with a Lebesgue probability density function $p_n(\tau^2),\ \tau^2> 0,$ that satisfies the following three conditions:
  \begin{itemize}
    \setlength\itemsep{0em}
   \item[a)] The density $p_n(\tau^2)$ is monotonically decreasing in $\tau^2$.
   \item[b)] There exists a constant $c_1>0$ such that $p_n(c_1(\tau^2)^\ast)\geq \exp\left(-n\epsilon_n^2\right)$ for all sufficiently large $n$ where $(\tau^2)^{\ast}=(\tau^2)^{\ast}(n)=n^{-1/(2m_0+1)}$.
   \item[c)] There exists a constant $c_2>0$ such that $\int_{c_2(\tau^2)^\ast}^\infty p_n(\tau^2)d\tau^2\leq \exp(-5\ n\epsilon_n^2)$ for all sufficiently large $n$ with $(\tau^2)^\ast$ as in b).
   \end{itemize}
   \item[A6)] The parameter $\tau_{poly}^2$ that controls the strength of the penalty for the polynomial part is equal to some positive constant, i.e. $\tau_{poly}^2>0$.
   \item[A7)] The design points $x_1,\dots,x_n\in[0,1]$ are distributed in a regular manner, i.e. there exists a fixed distribution function $Q$ that is associated with a positive and continuous design density $\rho$ such that $\|Q-Q_n\|_\infty=o(k_0^{-1})$, where $Q_n(x)=n^{-1}\sum_{i=1}^n \dsOne_{\{x_i\leq x\}}$, $x\in[0,1],$ is the empirical cumulative distribution function of the design points. 
 \end{itemize}

Similar assumptions as A1)--A4) and A7) have been used by \cite{ClaKriOps2009} and \cite{Xia2019} for their analyses of the frequentist O-splines estimator~\eqref{EstimatedCoeffVector}. Assumptions A5) and A6) are unique to the present Bayesian setting. Assumption A5) refers to the hyperprior on the smoothing variance and is particularly interesting. In Section~\ref{subsec:OptimalHypis} we investigate which of the hyperpriors $p(\tau^2)$ investigated in Proposition~\ref{MarginalPrior} satisfy Assumption A5). Before that, we state our main posterior concentration result for Bayesian penalized splines, Theorem~\ref{theo}, along with a sketch of our proof.

\subsection{Main theorem and sketch of the proof} 

\begin{theorem}\label{theo}
 Under {Assumptions} {A1)}--{A7)}, the O-splines posterior $\widetilde{\Pi}_n\left(\cdot\mid Y^n,\sigma_0^2\right)$ concentrates around $f_0$ at rate $\epsilon_n=n^{-m_0/(2m_0+1)} (\log n)^{1/2}$, i.e.~there exists a constant $M>0$ such that
\begin{align}\label{PosteriorConcentration}
   \dsE_0 {\widetilde{\Pi}}_n(\{f\in C([0,1]):\|f-f_0\|_n\geq M\epsilon_n\}\mid Y^n,\sigma_0^2)=o(1),
\end{align}
where $\dsE_0$ is the expectation under the true model \eqref{GaussiannonparametricRegressionModel}.
\end{theorem}

\subsection*{Sketch of the proof}

A detailed proof of Theorem~\ref{theo} is provided in Section~D of the Supplement. In the following, we sketch the main steps of this proof: As mentioned before, we use the KL and testing strategy for our proof of Theorem~\ref{theo}. Thus, we need to verify the KL condition and the testing condition. If both conditions are satisfied, then~\eqref{PosteriorConcentration} follows from general arguments \citep[see, e.g.,][]{Rou2016}. The proper O-splines prior~\eqref{ProperOSplinesPrior} enables us to apply the KL and testing strategy in the present setting. This would not be possible for the usual partially improper O-splines prior~\eqref{OSplinesPrior} as the KL and testing strategy requires a proper prior.

\vspace{-1em}
\paragraph{KL condition:} Intuitively, the KL condition ensures that the prior is sufficiently ``thick'' around the true function $f_0$. To verify the KL condition, we need to show that the prior probability of KL type neighbourhoods of $f_0$ is bounded from below, i.e.~we need to show that there exists a constant $c^\ast>0$ such that
      \begin{align}\label{PriorProbabilityKLTypeNeighbourhoods}
   \widetilde{\Pi}_n\left(\{f\in C([0,1]):K(f_0\|f)\leq n\epsilon_n^2,\, V(f_0\|f)\leq n\epsilon_n^2\}\right)\geq \exp(-c^\ast n\epsilon_n^2),
  \end{align}
  where $K(f_0\|f)$ is the KL divergence and $V(f_0\|f)$ is the KL variation between the true function $f_0$ and a continuous function $f\in C([0,1])$. The key idea of our proof is to exploit the simple conditionally Gaussian form~\eqref{ProperOSplinesPriorDemmlerReinsch} of the proper conditional O-splines prior in terms of the coefficients $u\in\dsR^d$ of a DR basis $\{Z_1,\dots,Z_d\}$. Omitting many technical details, we proceed as follows: 
  
  First, we introduce a sequence of DR bases $\{Z_{n,1},\dots,Z_{n,d}\}$ for the sequence of spline spaces $S(m,\xi_n)$. The existence of such a sequence of DR bases is established in Proposition~C.1 in the Supplement. This sequence of DR bases gives rise to a sequence of pseudo-true coefficient vectors $u_0=u_0(n),\ n\in\dsN,$ where $u_0=(Z^\top Z)^{-1}Z^\top f_0^n=Z^\top f_0^n/n$ is the coefficient vector of the empirical proximum of $f_0$ in $S(m,\xi_n)$. Using standard results from spline approximation theory \citep[][Theorem 20.3]{Pow1981}, we conclude that it suffices to bound the prior probability of a shrinking sequence of Euclidean balls around $u_0$ from below.
  
  Using~\eqref{ProperOSplinesPriorDemmlerReinsch} we split the prior on the DR coefficient vector $(u_1,\dots,u_d)$ into a polynomial part and an independent non-polynomial part. The polynomial part follows a Gaussian prior $N_q(0,\tau_{poly}^2I_q)$, which facilitates the derivation of a lower bound for the prior probability. To this end, we use a general inequality for the probability of a shifted Euclidean ball under a centered Gaussian measure that we have adapted from \cite{HofSheDud1979} (see Lemma~C.4 in the Supplement).
      
      For the non-polynomial part we apply the law of total probability, which gives us an integral over the hyperprior $p_n(\tau^2)$ on the smoothing variance. To obtain a lower bound for this integral, we use Assumptions A5)a) and A5)b) and then exploit the simple conditionally Gaussian form~\eqref{ProperOSplinesPriorDemmlerReinsch}. For the latter, we use again the general inequality adapted from \cite{HofSheDud1979}.
      
      Using Markov’s inequality and bounds for the triangular array of DR eigenvalues $\gamma^{(q)}_{n,j},\ j=1,\dots,d,$ from \cite{Xia2019}, we find that the prior probability of the KL type neighbourhoods~\eqref{PriorProbabilityKLTypeNeighbourhoods} is bounded from below by $\exp(-c^\ast n \epsilon_n^2)$ with $c^\ast=2$.
      
      \vspace{-1em}
    \paragraph{Testing condition:} Intuitively, the testing condition bounds the complexity of the effective support of the prior $\widetilde{\Pi}_n$. To verify the testing condition, we introduce the test
\begin{align*}
    \phi_n=\dsOne\{\|\widehat{f}_t-f_0\|_n\geq M/2\ \epsilon_n\},
\end{align*}
where $\widehat{f}_t=\widehat{f}_t(n)$ is a truncated DR estimator with cut-off ${t=q+\lceil n^{1/(2m_0+1)}\rceil}$ and $M>0$ is a positive constant. In addition to that, we introduce the
    subsets $\Theta_n=\lbrace f\in \mathcal{S}(m,\xi_n):\sum_{j=t+1}^d u_j^2 \leq \widetilde{M}\epsilon_n^2\rbrace,$ where $\widetilde{M}>0$ is another positive constants. The subsets $\Theta_n$ contain less wiggly splines in $\mathcal{S}(m,\xi_n)$ as the squared norm of the high-frequency DR coefficients $(u_{t+1},\dots,u_d)$ is bounded from above. Following \cite{Rou2016}, we need to show that:
       \begin{enumerate}[i)]
   \item The prior $\widetilde{\Pi}_n$ is effectively supported on $\Theta_n$, i.e. $\widetilde{\Pi}_n\left(\Theta_n^C\right)=o\left(\exp(-4\ n\epsilon_n^2)\right)$.
   \item The type I error probability of the test $\phi_n$ vanishes, i.e. $\dsP_0(\phi_n=1)=o(1)$. 
   \item Restricted to the effective support $\Theta_n$, the type II error probability of the test $\phi_n$ vanishes sufficiently fast, i.e. 
   \begin{align*}
    \underset{f\in\Theta_n:\lVert f-f_0\rVert_n\geq M\epsilon_n}{\text{sup}} \dsP_{f}(\phi_n=0)=o\left(\exp(-4\ n\epsilon_n^2)\right).
   \end{align*}
   \end{enumerate}
   
For i), we use the law of total probability and Assumption A5)c). Then we use~\eqref{ProperOSplinesPriorDemmlerReinsch} and a tail bound for a linear combination of iid $\chi^2(1)$-variates from \cite{LauMas2000}.

For ii) we first apply Markov’s inequality such that
\begin{align*}
  \dsP_0(\phi_n=1)=\dsP_0\left(\lVert \widehat{f}_t-f_0\rVert_n\geq {M}/{2}\ \epsilon_n\right)\leq\dfrac{\dsE_0\lVert \widehat{f}_t-f_0\rVert_n^2}{M^2/4\ \epsilon_n^2}.
\end{align*}
Hence, it suffices to show that the truncated DR estimator $\widehat{f}_t$ is able to achieve the optimal rate of convergence, i.e.~$\dsE_0\|\widehat{f}_t-f_0\|_n^2=O\left(n^{-2m_0/(2m_0+1)}\right).$
For the latter result we proceed similar to \cite{ClaKriOps2009} for their analysis of the frequentist O-splines estimator~\eqref{EstimatedCoeffVector} and decompose the average mean squared error $\dsE_0\|\widehat{f}_t-f_0\|_n^2$ in a variance term, a shrinkage bias term and an approximation bias term. The variance term and the approximation bias term are relatively easy to control. For the shrinkage bias term, we need to show that the high-frequency coefficients of the empirical spline proximum are bounded, i.e. we need to show that $\sum_{j=t+1}^d u_{0,j}^2=O\left(n^{-2m_0/(2m_0+1)}\right)$. To this end, we use the results of \cite{Xia2019} for the frequentist O-splines estimator~\eqref{EstimatedCoeffVector} in a clever way (see Proposition~C.2 in the Supplement for details).

For iii) we use Proposition~\ref{propType2err} and a tail bound for a $\chi^2(t)$-variate from \cite{ArmDunLee2013b}. As explained in Section~\ref{sec:tDRB}, controlling the type II error probability is much more difficult for the frequentist O-splines estimator~\eqref{EstimatedCoeffVector}, which explains why we needed to introduce the novel truncated DR estimator $\widehat{f}_t$ in this paper. 

The KL condition and the testing condition together imply~\eqref{PosteriorConcentration} which finishes the proof of Theorem~\ref{theo}. A more detailed proof is provided in Section~D of the Supplement. $\square$

\subsection{Hyperpriors leading to posterior concentration at near optimal rate}\label{subsec:OptimalHypis}
Under Assumption {A1)}, the contraction rate $\epsilon_n=n^{-m_0/(2m_0+1)}(\log n)^{1/2}$ in Theorem~\ref{theo} is \emph{near optimal} in the sense that it is the minimax rate \citep{Sto1982} times a logarithmic  factor. Logarithmic factors such as $(\log n)^{1/2}$ are common in the context of posterior concentration rates \cite[e.g.~][]{Rou2016}.

Theorem~\ref{theo} shows that we can achieve posterior concentration at near optimal rate for the Bayesian O-splines approach under certain conditions on the hyperprior $p(\tau^2)$. In the following Corollary~\ref{cor1} we investigate which of the hyperpriors $p(\tau^2)$ for the smoothing variance that were considered in Proposition~\ref{MarginalPrior} satisfy Assumption {A5)} and therefore lead to posterior concentration at near optimal rate.

\begin{corollary}\label{cor1}
 The following hyperpriors $p_n(\tau^2)$ satisfy Assumption {A5)}.
 \begin{enumerate}  \setlength\itemsep{1em}
 \item {Uniform}: $\tau^2\sim U(0,c(\tau^2)^\ast)$, where $c>0$ is some positive constant and $(\tau^2)^{\ast}=n^{-1/(2m_0+1)}$ as in Assumption {A5)}. 
  \item {Gamma}: $\tau^2\sim Ga(\alpha,\beta_n)$ with constant shape parameter $\alpha\in(0,1]$ and increasing rate parameter $\beta_n=c_{\beta}\ n\epsilon_n^2/(\tau^2)^\ast=c_{\beta}\ n^{2/(2m_0+1)}\log n$, where $c_{\beta}>0$ is some positive constant.
  \item {Weibull}: $\tau^2\sim \text{Weibull}(k,\lambda_n)$ with constant shape parameter $k\in(0,1]$ and increasing rate parameter $\lambda_n=c_\lambda\ (n\epsilon_n^2)^{1/k}/(\tau^2)^\ast$, where $c_\lambda>0$ is a positive constant.
 \end{enumerate}  \setlength\itemsep{1em}
 The following hyperpriors $p_n(\tau^2)$ do not satisfy Assumption A5).
 \begin{enumerate}
  \item {Inverse Gamma}: $\tau^2\sim \text{IG}(\alpha,\beta_n)$ with shape parameter $\alpha>0$ and scale parameter $\beta_n>0$.
  \item {Scaled Beta Prime}: $\tau^2\sim \mathit{SBP}(1,1,1/\lambda_n)$ with rate parameter $\lambda_n>0$.
  \end{enumerate}
\end{corollary}

\subsection{General remarks}

\begin{itemize}
\item The results of this section prove that we can achieve posterior concentration at near optimal rate for the Bayesian O-splines approach. Assumption A3) shows that we can use a faster rate for the number of spline knots $k_0$ than the slow regression spline rate $n^{1/(2m_0+1)}$ that is commonly being used for posterior concentration results in Bayesian spline regression \citep[e.g.][Section 7.7]{Ghovan2007}. Intuitively, the roughness penalty $\|D^qf\|^2_{\mathcal{L}^2}$ in the O-splines prior drags the posterior mass away from overly wiggly spline functions so that overfitting can be avoided.

\item From Assumption A5) it becomes clear that the hyperprior $p_n(\tau^2)$ on the smoothing variance $\tau^2$ needs to strike a fine balance to achieve posterior concentration at near optimal rate. It needs to be sufficiently thick at $c_1(\tau^2)^\ast$ to avoid underfitting (oversmoothing) but at the same time it needs to be sufficiently thin on $[c_2(\tau^2)^\ast,\infty)$ to avoid overfitting (undersmoothing).

\item Corollary~\ref{cor1} shows that the Uniform,  Gamma and  Weibull hyperpriors satisfy Assumption A5) and therefore lead to posterior concentration at near optimal rate according to Theorem~\ref{theo}. For those hyperpriors that do not satisfy Assumption {A5)}, we cannot make this conclusion. As the above conditions are only sufficient, we can, of course, not completely rule the latter hyperpriors out. Nevertheless, it seems reasonable to rather use a hyperprior for which optimal asymptotic performance is guaranteed. \cite{ArmDunLee2013b} have argued similarly for the screening of shrinkage priors in a high-dimensional linear regression context.

\item Our results cover in particular a Weibull hyperprior with shape parameter $1/2$. \cite{KleKne2016} compared various different hyperprior families $p(\tau^2)$ for the smoothing variance empirically and identified a Weibull hyperprior with shape parameter $1/2$ as optimal default, which they referred to as \emph{scale-dependent hyperprior}. Thus, the theoretical results of this section complement their recommendation based on extensive simulations.  
\end{itemize}

\subsection{Unknown residual variance}\label{SecUnknownVar}

Throughout, we have assumed that the true value $\sigma_0^2$ of the residual variance in model \eqref{GaussiannonparametricRegressionModel} is known. In this section, we explain how the Bayesian O-splines approach and our posterior concentration results can be generalized to cover the more realistic case when $\sigma_0^2$ is unknown. 

\begin{figure}[htbp]
  \begin{center}
\begin{tikzpicture}

  \node[obs](y) {$Y^n$};
  \node[latent, left=of y] (b) {${b}$};
  \node[latent, left=of b] (tau2){$\tau^2$};
  \node[latent, right=of y] (sigma2) {$\sigma^2$};
  
  \edge {tau2} {b} ; 
  \edge {b} {y} ; 
  \edge {sigma2} {y} ; 

\end{tikzpicture}
  \end{center}
  \caption{Graphical model of the joint O-splines prior in model \eqref{GaussiannonparametricRegressionModel}. Note that we do not condition the B-spline coefficients $b$ onto the residual variance $\sigma^2$ a priori, i.e. we do not use a prior of the form $b\mid\sigma^2$, but rather treat $b$ and $\sigma^2$ as independent a priori, which is recommended by \cite{MorRocGeo2019} to avoid biased estimation.}\label{figJointOSplinesPrior}
\end{figure}
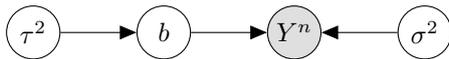

In the unknown $\sigma_0^2$ case, we endow the B-spline coefficients $b\in\dsR^d$ with the O-splines prior \eqref{OSplinesPrior} and the residual variance with an {independent} Inverse Gamma prior $\sigma^2\sim IG(\alpha_\epsilon,\beta_\epsilon)$. The fully Bayesian model specification then reads as follows
\begin{align}\label{FullHierarchy}
    Y^n\mid b,\sigma^2 &\sim N_n(Bb,\sigma^2I_n),\nonumber\\
         b\mid \tau^2&\sim p(b\mid \tau^2) \propto  (\tau^2)^{-(d-q)/2}\ \exp\left(-b^\top R_B^{(q)}b/(2\tau^2)\right),\nonumber\\
     \tau^2&\sim p(\tau^2),\nonumber\\
     \sigma^2&\sim IG(\alpha_\epsilon,\beta_\epsilon). 
\end{align}
Figure~\ref{figJointOSplinesPrior} shows a graphical model of the resulting \emph{joint O-splines prior} in model~\eqref{GaussiannonparametricRegressionModel}.

Posterior sampling is hardly more difficult when moving from known $\sigma_0^2$ to unknown $\sigma_0^2$: To this end, we simply include an additional step in the Gibbs sampler and draw in an alternating manner from the three full conditional posteriors $b\mid Y^n,\sigma^2,\tau^2$ and $\tau^2\mid Y^n,b,\sigma^2$ as well as $\sigma^2\mid Y^n,b,\tau^2$. The first two are exactly as before and the third full conditional posterior $\sigma^2\mid Y^n,b,\tau^2$ is simply an Inverse Gamma distribution. With much more effort, we can also generalize our posterior concentration results to the unknown $\sigma_0^2$ case. In the following we explain how this can be done, the details are omitted. To establish posterior concentration results in the unknown $\sigma_0^2$ case, we proceed as follows:

\begin{itemize}
    \item First, we need to define a suitable notion of \emph{joint posterior concentration}: We say that the joint posterior concentrates at rate $\epsilon_n$ around the truth $(f_0,\sigma_0^2)$, if the expected posterior probability of joint neighbourhoods of the form
 \begin{align}\label{JointPosteriorConcentration}
  \{(f,\sigma^2)\in C([0,1])\times (0,\infty):\lVert f-f_0\rVert_n< M\epsilon_n \land \left|\sigma^2/\sigma_0^2-1\right|< M\epsilon_n\}
 \end{align}
 converges to $1$. Joint posterior concentration in this sense or comparable has been studied by \cite{ChoSch2007,WeiReiHop2020,BaiMorAnt2020}. However, none of these papers has considered Bayesian penalized splines with a roughness penalty. 
 \item To prove joint posterior concentration in the sense of~\eqref{JointPosteriorConcentration} for the joint O-splines prior depicted in Figure~\ref{figJointOSplinesPrior}, we replace the partially improper O-splines prior \eqref{OSplinesPrior} on the B-spline coefficients $b\in\dsR^d$ by the proper version \eqref{ProperOSplinesPrior}. Thus, we consider the \emph{joint proper O-splines prior} on the pair $(b,\sigma^2)$ where $b$ follows the proper O-splines prior \eqref{ProperOSplinesPrior} and $\sigma^2$ follows an independent $IG(\alpha_\epsilon,\beta_\epsilon)$ prior. The propriety of this prior allows us to use the KL and testing strategy again.
\item The KL condition can be verified by introducing rectangles that are contained within the resulting KL type neighbourhoods. Furthermore, the test
     \begin{align*}
    \phi_n=\dsOne\left\lbrace\lVert \widehat{f}_t-f_0\rVert_n\geq {M}/{2}\ \epsilon_n\ \lor \left|{\widehat{\sigma^2_t}}/{\sigma_0^2}-1\right|\geq {M}/{2}\ \epsilon_n\right\rbrace,
   \end{align*}
   can be used for the testing condition. Thereby, $\widehat{f}_t$ is a truncated DR estimator as before and $\widehat{\sigma^2_t}={\lVert Y^n-\widehat{f}_t^n\rVert_2^2}/({n-t})$ is the corresponding estimator of the residual variance. 
\end{itemize}

 With this, we find that the joint O-splines posterior concentrates at near optimal rate around the truth $(f_0,\sigma_0^2)$ in the sense of \eqref{JointPosteriorConcentration} under almost the same assumptions as before if we use an independent Inverse Gamma prior $\sigma^2\sim IG(\alpha_\epsilon,\beta_\epsilon)$ with constant shape $\alpha_\epsilon\in(0,1]$ and constant scale $\beta_\epsilon>0$ on the residual variance $\sigma^2$. We only get the restriction $m_0\geq 2$ for the regularity of the unknown regression function $f_0\in C^{m_0}([0,1])$. To cover the case $m_0=1$, the prior on $\sigma^2$ needs to have even less prior mass close to the origin compared to an Inverse Gamma distribution. For instance, the square root of an Inverse Gamma distribution works well.

\section{Empirical Evidence}\label{sec:EmpEvidence}

In this section, we conduct a small simulation study to answer the following two questions:
\begin{enumerate}
\item What are the empirical benefits of using a hyperprior $p(\tau^2)$ on the smoothing variance as opposed to fixing the smoothing variance to some constant in advance?
\item Are there empirical arguments in favor of the proper O-splines prior \eqref{ConditionalProperOSplinesPrior} based on projections compared to the more established proper prior based on the popular mixed model reparametrization (MMR)?
\end{enumerate}

\subsection{Adaptivity of the hyperprior on the smoothing variance}\label{subsec:Adaptivity}
To answer the first question, we consider the following simulation experiment contrasting the Bayesian O-splines approach with a hyperprior $p(\tau^2)$ on the smoothing variance with different fixed values for the latter.

\vspace{-1.0em}
\paragraph{Simulation design}
  We assume $x_i\overset{iid}{\sim}U(0,1)$, $i=1,\ldots,n$, $n=100$, and  simulate $R=100$ replicate data sets  from model \eqref{GaussiannonparametricRegressionModel} with $\sigma_0^2=(1/4)^2$ for each of the seven test functions $f_k(x)=\sin\left((k-1)\pi\ x\right)$, $k=1,\ldots,7$, of increasing roughness $\|D^2f_k\|^2_{\mathcal{L}^2}=1/2\ (k-1)^4\pi^4$.
  
The test functions are estimated using cubic Bayesian O-splines  with $k_0=20$ equidistant interior knots in the interval $[0,1]$ and a roughness penalty of order $q=2$. Moreover, we use an Inverse Gamma prior $\sigma^2\sim IG(1/1000,1/1000)$ for the residual variance. For the smoothing variance $\tau^2$ we consider:
\begin{enumerate}[a)]
    \item a Weibull hyperprior $\tau^2\sim Weibull(1/2,1/500)$ vs.
    \item a fixed value $\tau^2\in\{0.5,5,50,500,5.000\}$.
\end{enumerate} 

Our choice for the Weibull hypeprior in a) is motivated as follows: \cite{KleKne2016} compared various  hyperpriors $p(\tau^2)$ empirically for Bayesian P-splines. The authors identified a Weibull hyperprior with shape parameter $1/2$ as optimal default, which they refer to as \emph{scale-dependent hyperprior}. Given the close relation of P-splines and O-splines \citep[see, e.g.,][]{Xia2019}, there is no need to repeat their analysis in the present setting so that we follow their recommendation. Note further that our posterior concentration results support the choice of the Weibull hyperprior (cf.~Corollary~\ref{cor1}). The choice $\lambda=1/500$ for the rate parameter  ensures that the scale of the induced prior on the splines $\mathcal{S}(m,\xi)$ roughly matches the scale of the test functions $f_k,\ k=1,\dots,7$, here ranging from $-1$ to $+1$. 

Posterior sampling can be done with standard Gibbs steps for $b$ and $\sigma^2$ from their full conditional posteriors. For $\tau^2$ in case a) we obtain a Generalized Exponential-Type full conditional posterior \citep{ProMab2010}, for which we use a Metropolis-Hastings step based on iteratively weigthed least square proposals as in \citet{KleKne2016}. We refer to Sections~A and~E of the Supplement for more information on the Generalized Exponential-Type distribution and posterior sampling, respectively.

We estimate the underlying test functions $f_k(x)=\sin\left((k-1)\pi\ x\right)$, $k=1,\ldots,7,$ using the posterior mean $\widehat{f}=\sum_{j=1}^d B_j\widehat{b_j}$ with $\widehat{b}=\dsE(b\mid Y^n)$ and compare the Weibull hyperprior and the different fixed values $\tau^2\in\{0.5,5,50,500,5.000\}$ in terms of their MSE$=\|\widehat{f^n}-f_k^n\|_2^2/n$.

\vspace{-1.0em}
\paragraph{Results}

\begin{figure}[htbp]
\begin{center}
    \includegraphics[width=\textwidth]{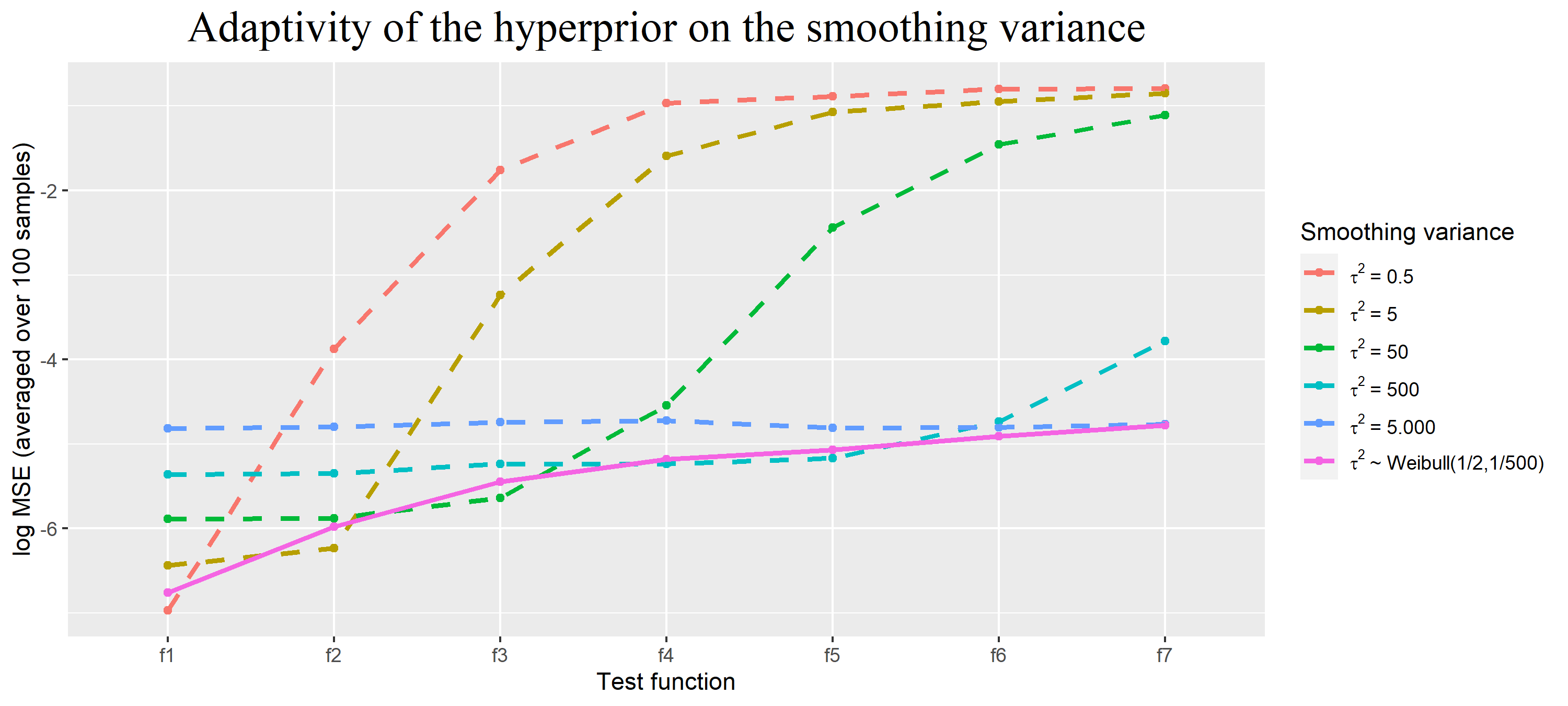}
\caption{Shown are the average log MSEs for Bayesian O-splines with a Weibull hyperprior $\tau^2\sim Weibull(1/2,1/500)$ on the smoothing variance (solid pink line) and for Bayesian O-splines with fixed values $\tau^2\in\{0.5,5,50,500,5.000\}$ for the latter. On the $x$-axis we go through the seven test functions $f_k(x)=\sin\left((k-1)\pi\ x\right)$, $k=1,\ldots,7$, of increasing roughness.}\label{figMSEhypi}
\end{center}
\end{figure}

Figure~\ref{figMSEhypi} shows the log MSE (averaged over the $100$ replicate data sets).
Figure~\ref{figMSEhypi} confirms our theoretical findings and the common practice in favor of using a hyperprior $p(\tau^2)$ on the smoothing variance as opposed to fixing the smoothing variance in advance. Figure~\ref{figMSEhypi} shows that a fixed smoothing variance is only slightly superior in terms of MSE for a very particular test function. Across the entire range of test functions, the results are much more robust when using a hyperprior. 
 
\subsection{Proper O-splines prior vs. proper prior based on MMR}\label{subsec:ProperPriors} 

To answer the second question, we investigate the proper O-splines prior \eqref{ConditionalProperOSplinesPrior} empirically and contrast it with the alternative definition based on the popular MMR. We will see that our proper prior \eqref{ConditionalProperOSplinesPrior} based on projections shows a much more reasonable behavior when varying the additional tuning parameter $\tau^2_{poly}$.

Recall that the usual conditional O-splines prior \eqref{ConditionalOSplinesPrior} is partially improper as the roughness penalty $\|D^qf\|^2_{\mathcal{L}^2}=b^\top R_B^{(q)}b$ does not penalize polynomials in $\mathcal{P}_{q-1}$. To establish posterior concentration rates, we have introduced the proper conditional O-splines prior \eqref{ConditionalProperOSplinesPrior}, which includes the additional penalty matrix $B^\top H^{(q)}B/(n\tau^2_{poly}),$ where $H^{(q)}=X_q(X_q^\top X_q)^{-1}X_q^\top$ is a projection matrix.

Previous definitions of a proper prior \citep[e.g.][]{WanOrm2008,HarRupWan2018,KleCarKneLanWag2021} are based on the MMR, which relies on an eigendecomposition of the roughness penalty matrix $R_B^{(q)}$ \citep[e.g.][]{FahKneLan2004b}. When using the MMR, the resulting proper conditional prior in terms of the B-spline coefficients $b\in\dsR^d$ is
\begin{align}\label{ProperPriorMMR}
   \Breve{p}(b\mid \tau^2)=N_d\left(b;0,\left({R_B^{(q)}}/{\tau^2}+Q_0Q_0^\top/{(n\tau^2_{poly})}\right)^{-1}\right),\ b\in\dsR^d,
 \end{align}
where $Q_0$ is a $d\times q$ matrix whose columns are orthonormal eigenvectors of $R_B^{(q)}$ corresponding to the $q$ zero eigenvalues of $R_B^{(q)}$. The proper prior~\eqref{ProperPriorMMR} can be regarded as an alternative to our proper prior~\eqref{ConditionalProperOSplinesPrior}. The prior~\eqref{ProperPriorMMR} can also be regarded as the Bayesian counterpart of the frequentist double penalty approach of \cite{MarWoo2011}. 
 
The proper conditional O-splines prior~\eqref{ConditionalProperOSplinesPrior} is advantageous for the derivation of posterior concentration rates because it leads to the simple form~\eqref{ProperOSplinesPriorDemmlerReinsch} in terms of the coefficients $u\in\dsR^d$ of a DR basis, which is not the case for the alternative proper prior~\eqref{ProperPriorMMR}. Next, we show that our proper prior~\eqref{ConditionalProperOSplinesPrior} is also advantageous from an empirical point of view. To this end, we consider the following simulation experiment.

\vspace{-1.0em}
\paragraph{Simulation design}
We simulate $n=100$ observations from model \eqref{GaussiannonparametricRegressionModel} with $\sigma^2_0=0.1^2$ and consider the two test functions $f_1(x)=4x-2$ (linear) and $f_2(x)=4x-2+\sin(2\pi x )$ (nonlinear) for a Uniform design $x_i\overset{iid}{\sim} U(0,1),\ i=1,\dots,n$. The test functions are estimated using cubic Bayesian O-splines with $k_0=20$ equidistant interior knots and a roughness penalty of order $q=2$. As in the previous Section~\ref{subsec:Adaptivity}, we use a Weibull hyperprior $\tau^2\sim Weibull(1/2,1/500)$ on the smoothing variance $\tau^2$ and an Inverse Gamma prior $\sigma^2\sim IG(1/1000,1/1000)$ on the residual variance $\sigma^2$. 

To ascertain the effect of the additional penalty matrices $B^\top H^{(q)}B/(n\tau^2_{poly})$ and $Q_0Q_0^\top/(n\tau^2_{poly})$ in the two proper priors~\eqref{ConditionalProperOSplinesPrior} and~\eqref{ProperPriorMMR}, we vary the parameter $\tau^2_{poly}$ and compute the resulting estimate $\widehat{f}=\widehat{f}(\tau^2_{poly})$ (posterior mean) under both proper priors~\eqref{ConditionalProperOSplinesPrior} and~\eqref{ProperPriorMMR}. To perform posterior estimation, we have adapted the MCMC sampler from the previous Section~\ref{subsec:Adaptivity} to cover the two proper priors~\eqref{ConditionalProperOSplinesPrior} and~\eqref{ProperPriorMMR} (see Supplement Section~E.2 for details).

\begin{figure}[htbp]
    \centering
    \includegraphics[width=\textwidth]{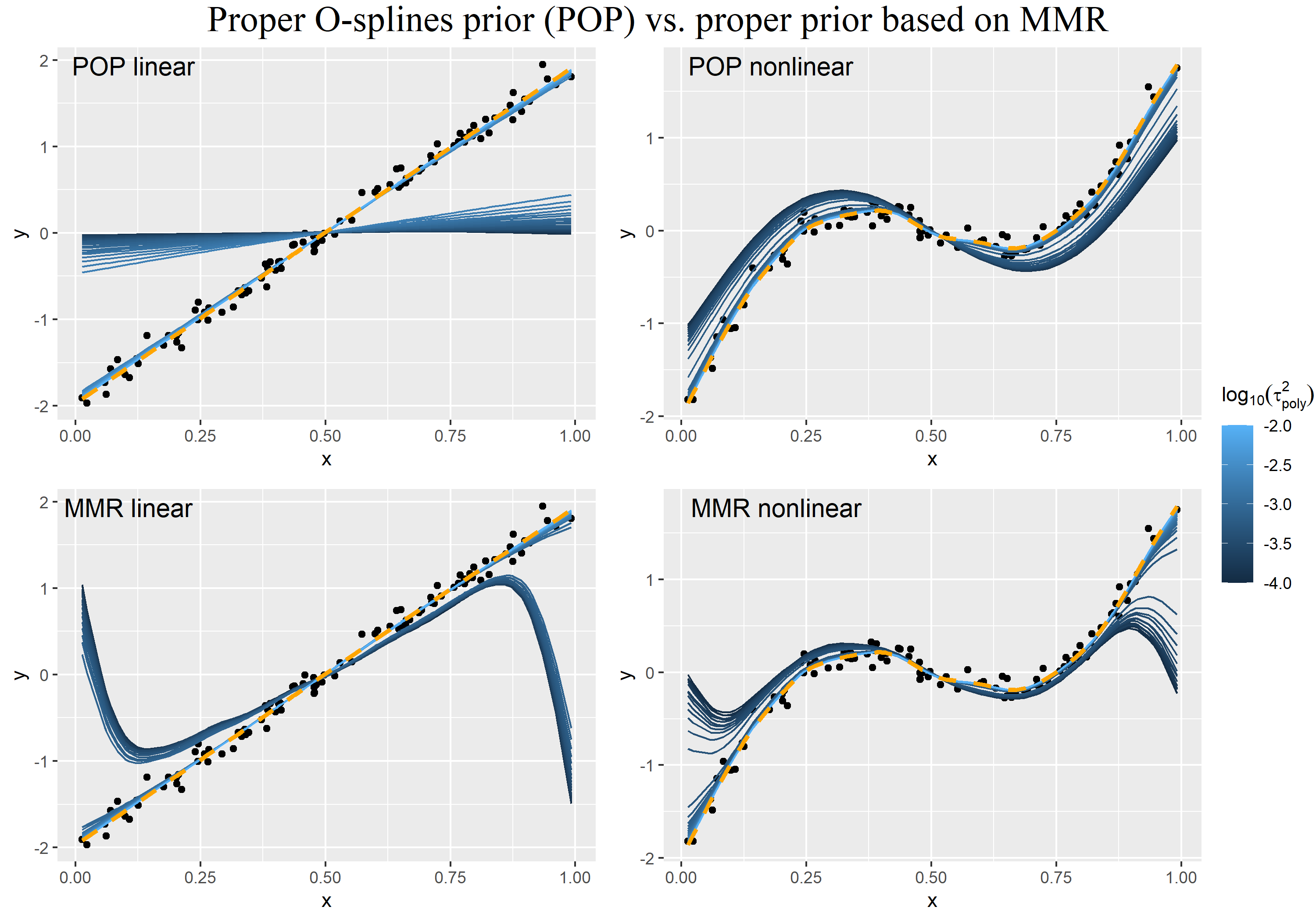}
    \caption{Shown are the estimates $\widehat{f}=\widehat{f}(\tau^2_{poly})$ for varying $\tau^2_{poly}$. The top panel shows the estimates for the proper O-splines prior (POP)~\eqref{ConditionalProperOSplinesPrior} based on projections, the bottom panel shows the estimates for the alternative proper prior~\eqref{ProperPriorMMR} based on the MMR for the linear test function $f_1(x)=4x-2$ (left column) and the nonlinear test function $f_2(x)=4x-2+\sin(2\pi x)$ (right column). The dashed orange line is the estimate $\widehat{f}$ obtained under the usual partially improper O-splines prior~\eqref{ConditionalOSplinesPrior}.}\label{fig:EstimatesProperPriors}
\end{figure}

\vspace{-1.0em}
\paragraph{Results}
    Figure~\ref{fig:EstimatesProperPriors} shows the estimates $\widehat{f}=\widehat{f}(\tau^2_{poly})$. Thereby, we vary the parameter $\tau^2_{poly}$ on the grid $\{10^{-4},10^{-3.95},10^{-3.9},10^{-3.85},\dots,10^{-2}\},$
    where we found a large sensitivity of $\widehat{f}$ to the precise value of $\tau^2_{poly}$. 
In the top panel of Figure~\ref{fig:EstimatesProperPriors}, we see that if $\tau^2_{poly}$ is not too small, i.e. $\tau^2_{poly}\geq 10^{-2.5}$, then the estimate obtained under the proper O-splines prior~\eqref{ConditionalProperOSplinesPrior} virtually coincides with the estimate obtained under the usual partially improper O-splines prior~\eqref{ConditionalOSplinesPrior} (dashed orange line). The same is true for the estimate obtained under the alternative proper prior~\eqref{ProperPriorMMR} based on the MMR, which can be seen in the bottom panel of Figure~\ref{fig:EstimatesProperPriors}.

However, as we decrease $\tau^2_{poly}$ and thereby increase the penalty for the linear part, the estimates obtained under the two proper priors~\eqref{ConditionalProperOSplinesPrior} and~\eqref{ProperPriorMMR} show a markedly different behavior: While the penalty matrix $B^\top H^{(q)}B/(n\tau^2_{poly})$ in the proper O-splines prior~\eqref{ConditionalProperOSplinesPrior} simply removes the linear trend from the estimate (top panel), the penalty matrix $Q_0Q_0^\top/(n\tau^2_{poly})$ in the alternative proper prior~\eqref{ProperPriorMMR} based on the MMR distorts the estimate in a rather strange manner (bottom panel). We want to emphasize that the observed behavior does not depend on the particular sample; it is largely independent of the underlying test function and it also occurs for a non-Uniform covariate design.

Overall, our results reveal a conceptual problem of the MMR: The penalized basis functions are not empirically orthogonal to the unpenalized basis functions. While the unpenalized basis functions span $\mathcal{P}_{q-1}$, the penalized basis functions do not span the orthogonal complement of $\mathcal{P}_{q-1}$ in $\mathcal{S}(m,\xi)$. The lack of orthogonality leads to undesirable consequences as depicted in Figure~\ref{fig:EstimatesProperPriors}, which has not always been recognized in the literature \citep[see, e.g.,][]{MarWoo2011,KleCarKneLanWag2021}.

\section{Discussion}\label{sec:Discuss}

Bayesian penalized splines are extremely popular in applied regression modelling but their theoretical understanding is lagging far behind. In this paper, we make an important step to close this gap in the literature and study posterior concentration rates for Bayesian O-splines in a Gaussian nonparametric regression model.

A key difference between Bayesian penalized splines and their frequentist counterpart is the additional hyperprior $p(\tau^2)$ on the smoothing variance $\tau^2$, which is only present in the Bayesian approach. Our results show that Bayesian O-splines can achieve posterior concentration at near optimal rate if the hyperprior strikes a fine balance between oversmoothing and undersmoothing, which can for instance be met by a Weibull hyperprior with shape parameter $1/2$. Following the reasoning of \cite{ArmDunLee2013b}, this is an argument in favor of such a hyperprior on the smoothing variance. The suitability of a Weibull hyperprior with shape parameter $1/2$ is also underlined by the empirical findings of \cite{KleKne2016} for Bayesian P-splines. However, to avoid oversmoothing the rate parameter of the Weibull hyperprior should in practice be set to a sufficiently small value or, even better, a prior scaling strategy as described in Section 2.4 of \cite{KleKne2016} should be used.

As a final note, we want to emphasize that our posterior concentration results can be generalized in many directions. First, it is very straightforward to establish similar results for the closely related Bayesian P-splines and Bayesian T-splines by combining our approach for Bayesian O-splines with the unified treatment of \cite{Xia2019} in the frequentist case. Second, our results can be carried over to multivariate or additive regression functions. Third, and most challenging, our insights may be combined with recent results of \citet{JeoGho2021,Bai2020b} to tackle non-Gaussian response models.
  Regarding the widespread use of Bayesian penalized splines in applied analyses, we pursue the generalization of our results in these directions in future research.


\bibliographystyle{dcu}
\renewcommand{\baselinestretch}{0.65}
\bibliography{main}  

\end{document}